\documentclass{article}

\usepackage{grffile}
\usepackage{ctable}

\usepackage{amssymb}
\usepackage{amsmath}
\usepackage{mathrsfs}
\usepackage{mathabx}
\usepackage[dvips]{epsfig}
\usepackage[small]{caption}
\usepackage{graphicx}
\usepackage[all]{xy}

\newtheorem{Lemma}{Lemma}[section]
\newtheorem{Theorem}{Theorem}
\newtheorem{Proposition}[Lemma]{Proposition}
\newtheorem{Corollary}[Lemma]{Corollary}
\newtheorem{Remark}[Lemma]{Remark}

\newtheorem{Hypothesis}[Lemma]{Hypothesis}

\newenvironment{Proof}%
 {\begin{trivlist} \item[]{\bf Proof. }}%
 {\hspace*{\fill}$\rule{.4\baselineskip}{.4\baselineskip}$\end{trivlist}}

 {\begin{trivlist}\item[]\textbf{Acknowledgments.}}{\end{trivlist}}

\makeatletter\@addtoreset{figure}{section}\makeatother

\makeatletter \@addtoreset{equation}{section} \makeatother

\newcommand{\R}{\mathbb{R}}
\newcommand{\C}{\mathbb{C}}
\newcommand{\N}{\mathbb{N}}
\newcommand{\Z}{\mathbb{Z}}

\def\Re{\mathop{\mathrm{Re}}}
\def\Im{\mathop{\mathrm{Im}}}

\newcommand{\rmO}{\mathrm{O}}

\newcommand{\rmd}{\mathrm{d}}
\newcommand{\rme}{\mathrm{e}}
\newcommand{\rmi}{\mathrm{i}}
\renewcommand{\ker}{\mathrm{Ker}\,}
\newcommand{\id}{\mathrm{\,id}\,}

\newcommand{\Rg}{\mathrm{Rg}}

\newcommand{\sm}{\mathscr{M}}

\newcommand{\cok}{\mathrm{Cok}\,}

\newcommand{\p}{\mathbb{P}}
\newcommand{\spn}{\mathrm{\,span}\,}

\makeatletter
\newsavebox{\@brx}
\newcommand{\llangle}[1][]{\savebox{\@brx}{\(\m@th{#1\langle}\)}%
  \mathopen{\copy\@brx\kern-0.5\wd\@brx\usebox{\@brx}}}
\newcommand{\rrangle}[1][]{\savebox{\@brx}{\(\m@th{#1\rangle}\)}%
  \mathclose{\copy\@brx\kern-0.5\wd\@brx\usebox{\@brx}}}
\makeatother

\newcommand{\up}{u_\mathrm{p}}
\newcommand{\upk}{u_{\mathrm{p},k_*}}

\renewcommand{\leq}{\leqslant}
\renewcommand{\geq}{\geqslant}

\newcommand{\rmnum}[1]{\romannumeral #1}
\newcommand{\Rmnum}[1]{\uppercase\expandafter{\romannumeral #1\relax}}

\def\Xint#1{\mathchoice
   {\XXint\displaystyle\textstyle{#1}}%
   {\XXint\textstyle\scriptstyle{#1}}%
   {\XXint\scriptstyle\scriptscriptstyle{#1}}%
   {\XXint\scriptscriptstyle\scriptscriptstyle{#1}}%
   \!\int}
\def\XXint#1#2#3{{\setbox0=\hbox{$#1{#2#3}{\int}$}
     \vcenter{\hbox{$#2#3$}}\kern-.5\wd0}}

\def\dashint{\Xint-}

\newfam\bifam
\font\tenbi=cmmib10 scaled \magstep1 \font\sevenbi=cmmib10 at 11pt
\font\fivebi=cmmib10 at 6pt \textfont\bifam = \tenbi
\scriptfont\bifam = \sevenbi \scriptscriptfont\bifam= \fivebi

\ifx\pdfoutput\undefined
   \pdffalse
\else
   \pdfoutput=1
   \pdftrue
\fi
\ifpdf
   \usepackage{graphicx}
   \usepackage{epstopdf}
\else
   \usepackage{graphicx}
\fi

\begin{document}

\begin{center}
{\fontsize{16}{16}\fontfamily{cmr}\fontseries{b}\selectfont{The Effect of Impurities on Striped Phases}}\\[0.2in]
Gabriela Jaramillo$\,^1$,  Arnd Scheel$\,^2$, and Qiliang Wu$\,^3$\\
\textit{\footnotesize $\,^1$The University of Arizona, Department of Mathematics, 617 N. Santa Rita Ave, Tucson, AZ 85721, USA\\
$\,^2$University of Minnesota, School of Mathematics,   206 Church St. S.E., Minneapolis, MN 55455, USA\\
$\,^3$Michigan State University, Department of Mathematics,  619 Red Cedar RD, East Lansing, MI 48824,  USA}
\date{\small \today} 
\end{center}

\begin{abstract}
\noindent 
We study the effect of algebraically localized impurities on striped phases in one space-dimension. We therefore develop a functional-analytic framework which allows us to cast the perturbation problem as a regular Fredholm problem despite the presence of essential spectrum, caused by the soft translational mode. Our results establish the selection of jumps in wavenumber and phase, depending on the location of the impurity and the average wavenumber in the system. We also show that, for select locations, the jump in the wavenumber vanishes. 
\end{abstract}

%

\vspace*{0.2in}

{\small
{\bf Running head:} {Impurities in  striped phases}

}
\vspace*{0.2in}

%
%
%
%

\section{Introduction}

We are interested in the effect of localized impurities on self-organized, spatially periodic patterns, in particular in the idealized situation of an unbounded domain. Our goal is to quantify the effect of the impurity on phases and wavenumbers in the far field. A prototypical example for the formation of self-organized periodic patterns is the Swift-Hohenberg equation
\[
u_t=-(\Delta +1)^2 u + \mu u - u^3,
\]
where, for $0<\mu\ll 1$, periodic patterns of the form $u_*(kx;k)$, $u_*(\xi;k)=u_*(\xi+2\pi;k)$, exist for a band of admissible wavenumbers $k\in (k_-(\mu),k_+(\mu))$. Our results are concerned with this system in one-dimensional space, $x\in\R$, including an impurity,
\begin{equation}\label{e:sh}
u_t=-(\partial_x^2 +1)^2 u + \mu u - u^3+\varepsilon g(x,u),
\end{equation}
where $|g(x,u)|\leq C(u)(1+|x|)^{-\gamma_*}$, for some $\gamma_*$ sufficiently large. 

We find such perturbation problems interesting for a variety of reasons. First, small impurities are simple examples of defects in spatially extended systems, and a systematic description of such defects is essential to various multi-scale descriptions of extended systems. In particular, defects can be responsible for the selection of wavenumbers $k$ in extended systems. Second, perturbations of periodic patterns pose challenging technical problems since the linearization at such periodic structures is generally not Fredholm when considered as an operator on translation-invariant (or algebraically weighted) function spaces. The difficulty stems from the presence of a non-localized neutral (or soft) mode, in this case the derivative $\partial_xu_*$ of the periodic pattern, which induces a branch of essential spectrum near the origin. In this regard, our results can be viewed as a continuation of a variety of results on perturbation and bifurcation in the presence of essential spectrum. Third, one can interpret the effect of inhomogeneities in relation to the notorious question of asymptotic stability of periodic patterns, where the pattern is perturbed at time $t=0$, whereas in our case the perturbation is constant in time. It would be quite interesting to bring those two view points together and study spatio-temporal perturbations of striped phases; see, for instance, \cite{dsss, gallay,zumbrunwith,johnsonzumbrun,uecker,scheelwu,schneider}. 

The effect of inhomogeneities on patterns with soft modes, that is, with eigenmodes of the linearization that exhibit neutral or weak temporal decay, has been studied in detail when periodic patterns are oscillatory in time \cite{kollar,SSdef}. In this case, inhomogeneities may create wave-sources such as target patterns, or act as weak sinks. In fact, in this case, the effects are quite similar to the effect of boundary conditions on oscillatory media, or, more generally, the effect of self-organized coherent structures on waves in the far-field. 

In the case of stationary periodic patterns, with vanishing group velocities, as they arise in the Swift-Hohenberg equation, the literature on defects and their characterization is quite extensive \cite{defectsSH}, albeit arguably not at the level of detail as we are striving for, here. In the direction of the present work, the characterization of boundary conditions on striped phases in \cite{morrissey} is closest. Results there show how to identify and compute strain-displacement relations, that is, relations between wavenumbers and phases (translations) of periodic patterns in the far field, induced by the presence of the boundary. Our present work can be viewed as matching such relations at $+\infty$ and $-\infty$.

Technically, our work is following up on recent studies of inhomogeneities in a variety of contexts \cite{jara3,jara1,jara2}, where Kondratiev spaces were used to study perturbations of spatio-temporally periodic patterns by inhomogeneities. The present work goes however significantly past those techniques by treating non-normal form, actual periodic patterns, where in \cite{jara3,jara1,jara2} the periodic patterns were, after appropriate transformations, constant in space. 

Our results are concerned with the spatially one-dimensional situation, only, but we hope that our approach will allow us to tackle higher-dimensional problems, as well. From a phenomenological point of view, the one-dimensional case is most difficult since effective diffusion of the neutral mode is weakest in one space-dimension, such that the effect of the inhomogeneity on the far-field is the most significant. This phenomenon is well understood in the case of diffusive stability, where decay of localized data is faster in $n$ space-dimensions $t^{-n/2}$, or in the case of impurities in oscillatory media, where small impurities can generate wave sources only in dimensions $n\leq 2$ \cite{jara3,jara2,kollar}. From a technical point of view, the one-dimensional case is easiest since the problem of finding stationary solutions can be cast as an ordinary differential equation; see for instance \cite{morrissey,SSdef} for this point of view. Our approach is different and in some sense more direct.  We will however comment on how to implement a proof using such ``spatial dynamics'' methods in the discussion. 

\paragraph{Notations}  We collect some useful notation. Let $\mathbb{P}_j(\R)$ and $\mathbb{P}_j(\Z)$ denote the set of complex-coefficient polynomials of degree less than $j\in\Z^+$ defined on the real line and on the set of integers, respectively. The inner product in a Hilbert space $H$ is denoted as $\langle\cdot, \cdot\rangle$ and the linear subspace spanned by $u\in H$ is denoted as $\langle u \rangle$. The Fourier transform on $L^2(\R, H)$ and $L^2(\Z, H)$ are denoted respectively as $\mathcal{F}$ and $\mathcal{F}_{\rm d}$. Moreover, for a Banach space $B$, the notation $\llangle u^*, u \rrangle$ represents the action of a linear functional $u^*\in B^*$ on $u\in B$. Throughout, the Lie bracket, $[L_1, L_2]$, of two operators $L_1$ and $L_2$ is the operator
\[
[L_1, L_2]:=L_1\circ L_2-L_2\circ L_1.
\]
We will use Banach spaces of functions on $\R$ and $\Z$.  Given $s\in\Z^+\cup\{0\}$, $p \in (1, \infty)$, $\gamma\in\R$, and denoting $\lfloor x \rfloor = \sqrt{ 1 +|x|^2}$, the weighted Sobolev space $W^{s,p}_\gamma$ is defined as
\[
W^{s,p}_\gamma:=\left\{u \in L^1_{\mathrm{loc}}(\R, H)\middle| \lfloor x \rfloor ^{\gamma}\partial_x^\alpha u \in L^p(\R, H), 
\text{for all }\alpha\in[0, s]\cap\Z\right\},
\]
with norm $\sum_{\alpha=0}^s\|\lfloor x \rfloor^{\gamma}\partial_x^\alpha u \|_{L^p}$, while the Kondratiev space $M^{s,p}_\gamma$ on $\R$ is defined as 
\[
M^{s,p}_\gamma:=\left\{u \in L^1_{\mathrm{loc}}(\R, H)\middle| \lfloor x \rfloor^{\gamma+\alpha}\partial_x^\alpha u \in L^p(\R, H), 
\text{for all }\alpha\in[0, s]\cap\Z\right\},
\]
 with norm $\sum_{\alpha=0}^s\| \lfloor x \rfloor^{\gamma+\alpha}\partial_x^\alpha u \|_{L^p}$. Their dual spaces are defined in the standard way and we write
 \[
 W^{-s,q}_{-\gamma}:=(W^{s,p}_\gamma)^*,\quad M^{-s,q}_{-\gamma}:=(M^{s,p}_\gamma)^*, \text{ where }1/p+1/q=1.
 \]
 For $s=0$, both spaces are simply weighted $L^p$-space, denoted as $L^p_\gamma$. For $p=2$, we denote $W^{s,2}_\gamma$ as $H^s_\gamma$.
 Additionally,  one can allow different weights on $\R^\pm$ to obtain an anisotropic version of these spaces. More specifically, letting $\chi_\pm$ be a smooth partition of unity, with $\mathrm{supp}(\chi_+)\subset (-1,\infty)$, $\chi_-(x)=\chi_+(-x)$, we define
 \[
 W^{s,p}_{\gamma_-, \gamma_+}:=\left\{u \in L^1_{\mathrm{loc}}(\R, H)\middle| \chi_\pm u\in W^{s,p}_{\gamma_\pm}\right\}, \quad
 M^{s,p}_{\gamma_-, \gamma_+}:=\left\{u \in L^1_{\mathrm{loc}}(\R, H)\middle| \chi_\pm u\in M^{s,p}_{\gamma_\pm}\right\},
 \] 
 which are Banach spaces respectively with norms 
 \[
 \|u\|_{W^{s,p}_{\gamma_-, \gamma_+}}:=\|\chi_+u\|_{W^{s,p}_{\gamma_+}} + \|\chi_-u\|_{W^{s,p}_{\gamma_-}},\quad
 \|u\|_{M^{s,p}_{\gamma_-, \gamma_+}}:=\|\chi_+u\|_{M^{s,p}_{\gamma_+}} + \|\chi_-u\|_{M^{s,p}_{\gamma_-}}.
 \]
Replacing $\R$ with $\Z$ and $\partial_x$ with the discrete derivative $\delta_+(\{u_j\}_{j\in\Z}):=\{u_{j+1}-u_j\}_{j\in\Z}$,  the discrete counterparts of $L^p_{\gamma_-,\gamma_+}$ and $M^{s,p}_{\gamma_-, \gamma_+}$ are denoted respectively as $\ell^p_{\gamma_-,\gamma_+}$, and $\sm^{s,p}_{\gamma_-, \gamma_+}$. We point out that the discrete counterparts of $W^{s,p}_{\gamma_-,\gamma_+}$ are isomorphic to $\ell^p_{\gamma_-,\gamma_+}$ due to the fact that $\delta_+$ is a bounded linear operator on $\ell^p_{\gamma_-,\gamma_+}$.
 
\paragraph{Outline.} The remainder of the paper is organized as follows. In Section \ref{s:2}, we present our main results. Section \ref{s:3} establishes Fredholm properties of one-dimensional differential operators with periodic coefficients in suitable algebraically weighted spaces. Section \ref{s:4} exploits these weighted spaces to treat impurities via an implicit function theorem and establishes expansions for solutions. We conclude with a discussion in Section \ref{s:5}.

\paragraph{Acknowledgment.} The authors acknowledge partial support through the National Science Foundation through  grants NSF-DMS-1311740 (AS) and NSF DMS-1503115 (GJ).

\section{Main Result}\label{s:2}

We state assumptions and main results.

\begin{Hypothesis}[Localization of impurity]\label{h:0}
We consider \eqref{e:sh} with smooth inhomogeneity $g(x,u)$ that is algebraically localized,
\begin{equation}\label{:gloc}
|\partial^{j_1}_x\partial^{j_2}_ug(x,u)|\leq (1+|x|)^{-\gamma_*}, j_1+j_2\leq 3,
\end{equation}
where $\gamma_*>6$.
\end{Hypothesis}
We next assume the existence of a periodic pattern.
\begin{Hypothesis}[Existence of stripes]\label{h:1}
We assume that there exists an even, periodic solution $\up$ with wavenumber $k_*>0$, $\up(\xi;k_*)=\up(\xi+2\pi;k_*)=\up(-\xi;k_*)$, to
\begin{equation}\label{e:ssh}
-(k_*^2\partial_\xi^2 +1)^2 u + \mu u - u^3=0,
\end{equation}
for some $\mu>0$, fixed.
\end{Hypothesis}
Note that this assumption is satisfied for $0<\mu\ll 1$, $|k_*-1|\ll 1$. 

The next assumption requires in particular that $\up$ is Eckhaus-stable. In order to state this assumption precisely, we introduce the family of Bloch-wave operators
\begin{equation}\label{e:bl}
L_\mathrm{B}(\sigma):=-\left(1+(\partial_x+\rmi\sigma)^2\right)^2+\mu-3\up^2(x), \quad \sigma\in [0,k_*),
\end{equation}
defined on $\mathcal{D}(L_\mathrm{B}(\sigma))=H^4_\mathrm{per}(0,2\pi/k_*)\subset L^2_\mathrm{per}(0,2\pi/k_*)$. Note that  all $L_\mathrm{B}(\sigma)$ have compact resolvent and depend analytically on $\sigma$ as closed operators with Fredholm index 0.

\begin{Hypothesis}[Stability of stripes]\label{h:2}
We assume that the periodic solution $\up$ is spectrally stable, that is, $0\in \mathrm{spec}(L_\mathrm{B}(\sigma))$ precisely for $\sigma=0$, when the eigenvalue $\lambda=0$ is algebraically simple, with eigenfunction $\up'$. For $\sigma \sim 0$, the expansion of the zero eigenvalue in $\sigma$ does not vanish at second order, $\lambda(\sigma)=\lambda_2\sigma^2+\rmO(\sigma^3)$, for some $\lambda_2\neq 0$.
 \end{Hypothesis}
We note that for $\mu\ll 1$, Eckhaus-stable patterns satisfy this hypothesis with $\lambda_2<0$ \cite{mielke}, and 
Eckhaus-unstable patterns do not, due to a kernel of $L_B(\sigma)$ for some $\sigma\neq 0$. On the other hand, long-wavelength unstable patterns may satisfy this assumption with $\lambda_2>0$; see for instance \cite{sslong}. We will give an expression for $\lambda_2$ in \eqref{e:d}.

\begin{Lemma}[Family of stripes]\label{l:family}
There exists a smooth family of stripe solutions, $\up(kx-\varphi;k)$, to \eqref{e:sh}, parameterized by wavenumber $k\sim k_*$ and phase $\varphi\in \R/2\pi\Z$.
\end{Lemma}
\begin{Proof}
We solve 
\[
-(1+k^2\partial_\xi^2)^2 u + \mu u - u^3=0,
\]
as an equation $H^4_\mathrm{per,even}\to L^2_\mathrm{even}$ using the implicit function theorem near $\up(\xi;k_*)$. The assumption that the kernel of $L_\mathrm{B}(0)$ is simple, spanned by $\up'$, odd, guarantees invertibility of the linearization. 
\end{Proof}
Our main result is as follows. 
\begin{Theorem}\label{t:1}
Assume Hypotheses \ref{h:0}--\ref{h:2}. Then there exists $\varepsilon_0$ and  a two-parameter family of stationary solutions to \eqref{e:sh} of the form 
\[
u(x;\varepsilon)=\sum_\pm \chi_\pm(x)\up((k_*+k_0\pm k_1)x-\varphi_0\mp\varphi_1; k_*+k_0\pm k_1)+w(x),
\]
where $w\in H^4_{\gamma_*},$ $\gamma_*>6$, and $\varphi_1,k_1$ are $C^1$-functions of $\varepsilon,k_0\in (-\varepsilon_0,\varepsilon_0)$, $\varphi_0\in \R$. Moreover,  $k_1$ and $\varphi_1$ have the leading-order expansions 
\begin{align}\label{e:k1}
k_1&=M_k(\varphi_0,0)\varepsilon+\rmO(\varepsilon^2),\\
\varphi_1&=M_\varphi(\varphi_0,0)\varepsilon+\rmO(\varepsilon^2),
\end{align}
where for the case $k_0=0$,
\begin{align}
M_k(\varphi_0,0)&=\frac{\pi \displaystyle \int_\R g(x, \up(k_*x-\varphi_0; k_*)) \cdot\partial_\xi \up(k_*x-\varphi_0; k_*)\,\rmd x  }{\lambda_2 k_*\int_0^{2\pi/k_*} (\partial_\xi \up(k_*x; k_*))^2 \rmd x},\\
M_\varphi(\varphi_0,0)&=\frac{\pi  \displaystyle \int_\R  g(x, \up(k_*x-\varphi_0; k_*)) \cdot [(x-\varphi_0/k_*)\partial_\xi \up(k_*x-\varphi_0;k_*)+\partial_k\up(k_*x-\varphi_0; k_*)]\,\rmd x }{\lambda_2k_* \int_0^{2\pi/k_*} (\partial_\xi \up(k_*x; k_*))^2 \rmd x}.
\end{align}

\end{Theorem}
We note that when the inhomogeneity is a gradient field, i.e. $g = \partial_uG(x,u)$, then  
\[
\dashint M_k\,\rmd \varphi_0:=\frac{1}{2\pi}\int_0^{2\pi}M_k(\varphi_0,0)\,\rmd \varphi_0=0,
\] 
and $M_k$ necessarily vanishes for certain relative phase shifts $\varphi_0$. We can therefore find relative phase shifts for which $k_1=0$. 
\begin{Corollary}
Assume that $g \in H^1_{\gamma_*} $, $\gamma_*>6$, $M_k(\varphi_*,0)=0$, and $M_k'(\varphi_*,0)\neq 0$. Then there exists $\bar{\varepsilon}, \bar{k}_0>0$ and a function $\phi_0(\varepsilon,k_0): [0,\bar{\varepsilon}]\times[0,\bar{k}_0]\to \R$ with $\phi_0(0,0)=\varphi_*$ such that the wavenumber difference $k_1$ from Theorem \ref{t:1} vanishes for $\varphi_0=\phi_0(\varepsilon,k_0)$. 
\end{Corollary}
\begin{Proof}
Scaling the equation \eqref{e:k1} by $\varepsilon$ we may write $k_1 = \varepsilon \bar{k}$ where 
\[ \bar{k}(\varepsilon; \varphi_0, k_0)= M_k(\varphi_0,k_0) + O(\varepsilon).\]
Our assumptions $M_k(\varphi_*,0)=0$, $M_k'(\varphi_*,0)\neq 0$ imply that $\bar{k}=0$ satisfies the conditions for the implicit function theorem, guaranteeing the results of the corollary.
The conditions on $g$ allow us to obtain a well defined value for $M'_k(\varphi,0)$ .
\end{Proof}

\section{Fredholm properties in weighted spaces near the essential spectrum}\label{s:3}
The results in this section can be viewed independently of the remainder of the paper. The difficulty of perturbing a striped pattern is due to the fact that the linearization is not Fredholm, which in turn can be attributed to the presence of essential spectrum at the origin, which in turn is induced by the non-localized eigenfunction $\up'$. It is well known that the linearization ``behaves'' in many ways like an effective diffusion. We therefore expect that the linearization at a periodic pattern possesses properties similar to the Laplacian $\partial_{xx}$. The Laplacian, on the other hand, while not Fredholm when posed as a closed, densely defined operator mapping $\mathcal{D}(\partial_{xx})\subset L^2\to L^2$, is Fredholm when posed as a closed, densely defined operator mapping $\mathcal{D}(\partial_{xx})\subset L^2_{\gamma-2}\to L^2_\gamma$, for $\gamma\not\in \{\frac{1}{2},\frac{3}{2}\}$. The goal of this section is to generally describe Fredholm properties of operators with translation symmetry in $\R$ or $\Z$ near points of the essential spectrum. The main restrictions are to one unbounded spatial direction, to ``algebraically simple'' points of the essential spectrum, and to non-critical weights $\gamma$. 
Throughout, we consider bounded operators, only. We will point out how these results imply Fredholm properties for more general operators. 

The outline for this section is as follows. We first consider operators with unbounded variable $x\in\R$ in Section \ref{s:3.1}, then show how to adapt in a straight-forward fashion to operators with unbounded direction $\ell\in\Z$ in Section \ref{s:3.2}. We finally show how to relate those results to Floquet-Bloch theory for operators on $x\in\R$ with periodic coefficients. We establish Fredholm properties for those operators in Section \ref{s:3.3}. For convenience, we recall Fredholm properties of $\partial_{xx}$ and of its discrete analogue in the appendix. 

%
%
%

\subsection{Operators with continuous translation symmetry}\label{s:3.1}

\paragraph{Setup --- operator symbols and essential spectrum.}
We consider bounded operators $\mathcal{L}$ on $L^2(\R,Y)$, where $Y$ is a complex separable Hilbert space, that possess a translation symmetry, that is, they commute with the action of translations on $L^2(\R,Y)$. The Fourier transform is an isomorphism of $L^2(\R,Y)$, and, due to translation symmetry, the induced operator $\hat{\mathcal{L}}$ on the Fourier space is a direct integral of multiplication operators with Fourier symbol $\hat{\mathcal{L}}=\int_{k\in\R}L(k)\rmd k$, that is,
\begin{equation}
\begin{matrix}
\hat{\mathcal{L}}: & \mathcal{D}(\hat{\mathcal{L}})\subset L^2(\R,Y) & \longrightarrow & L^2(\R,Y) \\
& u(k) & \longmapsto & L(k)u(k),
\end{matrix}
\end{equation}
with $L(k)$ linear and bounded on $Y$ for all $k\in\R$, see \cite{amann1
}.
Formally, we have $\mathcal{L}=L(-\rmi\partial_x)$. We denote the Banach space of bounded operators on $Y$ by $B(Y)$, 
\begin{Hypothesis}[Analyticity of symbol]\label{h:3.1}
We assume that $L(k)$ is analytic, uniformly bounded, with values in $B(Y)$,  in a strip $k\in \Omega_0:=\R\times (-\rmi k_\rmi,\rmi k_\rmi)$ for some $k_\rmi>0$. 
Moreover, we require that $L(k)$ is Fredholm for all $k\in\R$ and invertible with uniform bounds for $|\Re k|\geq k_0>0$ for some $k_0$ sufficiently large. 
\end{Hypothesis}
We mainly think of $L(k)$ rational, $L(k)=P(k)Q(k)^{-1}$, with matrix-valued polynomials $P$ and $Q$, where the values of $k$ such that $Q(k)$ is singular lie off the real axis. On the other hand, our results allow to include convolution operators with exponentially localized kernels. Specific examples are $\partial_{xx}(1-\partial_{xx})^{-1}$, $\partial_x(1+\partial_x)^{-1}$, $(-\mathrm{id}+K*)$, $K$ an exponentially localized kernel, or $(1+\partial_x^2)^2 (1-\partial_x^2)^{-2}$. 

Note that the spectrum of $\mathcal{L}$ is bounded, given through
\[
\mathrm{spec}_{L^2(\R,Y)}\mathcal{L}=\{\lambda \mid L(k)-\lambda \text{ not bounded invertible for some }k\in\R\}.
\]
In the case $Y=\R^n$, this can be more explicitly characterized through
\[
\mathrm{spec}_{L^2(\R,\R^n)}\mathcal{L}=\{\lambda\mid\mathrm{det}\,(L(k)-\lambda)=0\}.
\]
Since $L(k)$ is invertible for large $k$ and Fredholm for all $k\in\R$, $L(k)$ is Fredholm of index $0$ for all $k\in\R$ and the set of $k\in\R$ where $L(k)$ is not invertible is discrete.

We are interested in the case where $\mathcal{L}$ is not invertible.

\begin{Hypothesis}[Simple kernel]\label{h:3.2}
There exists a unique $k_*$ and a unique  (up to scalar multiples) $e_0\neq 0$ such that $L(k_*)e_0=0$. We then scale $\langle e_0, e_0 \rangle=1$. 
\end{Hypothesis}
In particular, $\lambda=0$ belongs to the essential spectrum of $\mathcal{L}$. 
We can assume without loss of generality that $k_*=0$, possibly conjugating $\mathcal{L}$ with the multiplication operator $\rme^{\rmi k_* x}$. We write $e_0^*$ for the kernel of the adjoint $L^*(0)$ with $\langle e_0^*, e_0^*\rangle=1$. 


\paragraph{Spatial multiplicities in the essential spectrum.}

We are interested in the unfolding of the zero-eigenvalue at $k=0$ for the family $L(k)$. We therefore view $L(k)$ as an analytic operator pencil and define the \emph{spatial multiplicity} as the multiplicity of $k=0$ as an eigenvalue of the operator pencil. Since such constructions are possibly not widely known, and the use here is less standard, we include the relevant constructions. 

Recall that, according to Hypothesis \ref{h:3.2}, the kernel of $L(0)$ is one-dimensional. 
\begin{Lemma}\label{c:eigenexp}
There exists $m>0$, maximal, and $e(k)=\sum_{j=0}^{m}e_jk^j$ such that
\begin{equation}\label{e:LE2}
L(k)e(k)=\lambda_mk^me_0^*+\rmO(k^{m+1}),
\end{equation}
or, equivalently,
\[
\sum_{j=0}^k L_je_{k-j}=0,\quad k=0,\ldots,m-1; \qquad   \lambda_m:=\left\langle \sum_{j=0}^{m-1} L_{m-j}e_j,e_0^*\right\rangle\neq  0.
\]
Here, we expanded  $L(k)=\sum_{j=0}^{m}L_j k^j+\rmO(k^{m+1})$. We refer to $m$ as the \emph{spatial multiplicity} of $\lambda=0$. 
\end{Lemma}
\begin{Proof}
Write $Q_0$ for the orthogonal projection onto $\mathrm{span} \{e_0^*\}$. We solve $L(k)(e_0+v)=z$ by decomposing 
\begin{align}
\langle L(k)(e_0+v),e_0^*\rangle&=z_1\label{e:gLS1}\\
(\id-Q_0)L(k)(e_0+v)&=z_2,\label{e:gLSeig}
\end{align}
where $z=z_1e_0^*+z_2$, $z_1\in\R$ and $z_2\in\Rg(\id-Q_0)$. Since $L(0)$ is Fredholm of index 0, $L(0):e_0^\perp\to (e_0^*)^\perp$ is an isomorphism, and the second equation \eqref{e:gLSeig}  can be solved using the implicit function theorem, with solution $v=v_*(k, z_2)$, where $|k|, |z_2|$ small. We then plug $v_*(k, z_2)$ into \eqref{e:gLS1}, yielding
\[
f(k, z_1, z_2):=\langle L(k)(e_0+v_*(k,z_2)),e_0^*\rangle-z_1=0.
\]
Due to the fact that $L(k)$ is invertible for all $k\neq 0\in\Omega_0$, the reduced analytic function $f(k,0,0)$ has non-trivial Taylor jet, that is, there exists $m\in\Z^+$ and $\lambda_m\neq0\in\C$ so that $f(k,0,0)=\lambda_mk^m+\rmO(k^{m+1})$. Taking $v=v_*(k,0)$, we have
\[
L(k)(e_0+v_*(k,0))=f(k,0,0)e^*_0=\lambda_mk^m e^*_0+\rmO(k^{m+1}).
\]
Letting $e(k)$ be the Taylor expansion up to order $\rmO(k^m)$ of $e_0+v_*(k,0)$, the claims follow quickly.
\end{Proof}
\begin{Remark}\label{r:mult}
In the case where $\lambda$ is an algebraically simple eigenvalue of $L(0)$, one can slightly modify the construction in the proof of Lemma \ref{c:eigenexp} and solve $L(k)e(k)=\lambda(k)e(k)$ together with $\langle e(k)-e_0,e_0\rangle =0$ using Lyapunov-Schmidt reduction in much the same way. The linearization with respect to $(e,\lambda)$ is onto and one finds the function $\lambda(k)$ which is of course the expansion of the ``temporal eigenvalue '' $\lambda$ in the Fourier parameter $k$. From this construction, one finds $\lambda(k)=\tilde{\lambda}_mk^m+\rmO(k^{m+1})$, for some $\tilde{\lambda}_m\neq 0$, with $m$ as in Lemma \ref{c:eigenexp}.
\end{Remark}

Since expansions typically do not converge globally, we introduce localized expansions as follows. Define the pseudo-derivative symbols
\begin{align}
D(k)&=\rmi k(1+\rmi k)^{-1},\nonumber\\
D_{C,m}(k)&= k\left(1+C\rmi k^m\right)^{-1},\label{e:pd} 
\end{align}
with associated operators $D(-\rmi\partial_x), D_{C,m}(-\rmi\partial_x)$. 
Here $C>0$ will eventually be chosen sufficiently large so that the norm of the bounded multiplier $D_{C,m}$ is arbitrarily small. Restricting to the strip
\[
\Omega_0(C,m):=\{k\in\Omega_0\mid |\Im k|\leq k_1:= \frac{1}{\sqrt[m]{2C}}\sin(\frac{\pi}{2m})\},
\] 
$D_{C,m}(k)$ is in fact analytic and uniformly bounded, that is, there exists a constant $C(m)$ such that
\[
\|D_{C,m}(k)\|\leq \frac{C(m)}{\sqrt[m]{C}},\quad  \text{ for all }k\in\Omega_0(C,m).
\]
\begin{Remark}
On the enlarged strip, $\{k\in\C\mid |\Im k|< \frac{1}{\sqrt[m]{C}}\sin(\frac{\pi}{2m})\}$,
the pseudo-derivative $D_{C,m}$ is analytic but not bounded. To obtain boundedness, we can restrict ourselves to any narrower strip, $\{k\in\C\mid |\Im k|< \frac{1}{\sqrt[m]{N C}}\sin(\frac{\pi}{2m})\}$, for any $N>1$. For convenience, we simply chose $N=2$ and $\Omega_0(C,m)\subset \Omega_0$, where the strip $\Omega_0$ is introduced in Hypothesis \ref{h:3.1}.
\end{Remark}
Note that replacing $k$ by $D_{C,m}(k)$ in the expansion of $e(k)$ does not alter its Taylor expansion up to order $m$. We therefore may define, for all $k\in\Omega_0(C,m)$,
\[
\tilde{e}(k):=\sum_{j=0}^{m} \left[D_{C,m}(k)\right]^je_j,
\]
such that
\begin{equation}\label{e:expmod}
L(k)\tilde{e}(k)=
\lambda_me_0^*k^m+\rmO(k^{m+1}).
\end{equation}
Repeating these considerations for the adjoint, we also find $e^*(k)=\sum_{j=0}^{m}e_j^*\bar{k}^j$ and define 
\[
\tilde{e}^*(k):=\sum_{j=0}^{m} \left[\overline{D_{C,m}(k)}\right]^je_j^*,
\]
so that 
\begin{equation}\label{e:expmodad}
L^*(k)\tilde{e}^*(k)=
\bar{\lambda}_me_0 k^m+\rmO(k^{m+1}).
\end{equation}
Since $L^*(k)$ is anti-analytic, $e^*(k)$ is anti-analytic, and we use the complex conjugate $\overline{D_{C,m}(k)}$ to guarantee that  $\tilde{e}^*(k)$ is anti-analytic.

\paragraph{Fredholm properties of $\mathcal{L}$.}
The main results on Fredholm properties of $\mathcal{L}$ are stated in the following theorem. 
\begin{Proposition}[Fredholm properties of $\mathcal{L}$]\label{p:f1}
Suppose the operator $\mathcal{L}$
satisfies Hypothesis \ref{h:3.1} and \ref{h:3.2}, with $k^*=0$. Let $m$ be the spatial multiplicity according to Lemma \ref{c:eigenexp}. Then, for  $\gamma_-,\gamma_+ \not\in \{1/2,3/2, \cdots, m-1/2\}$, the operator
\begin{equation}\label{e:superL}
\mathcal{L}:\mathcal{D}(\mathcal{L})\subset L^2_{\gamma_--m,\gamma_+-m}(\R,Y)\to  L^2_{\gamma_-,\gamma_+}(\R,Y),
\end{equation}
 is closed, densely defined, and Fredholm.
Moreover, setting $\gamma_{\max} = \max \{ \gamma_-, \gamma_+\} , \gamma_{\min} = \min\{ \gamma_-, \gamma_+\}$, we have that 
\begin{itemize}
\item for $\gamma_{\min} \in I_m:=( m-1/2,\infty)$, the operator \eqref{e:superL} is one-to-one with cokernel
\[ \cok(\mathcal{L}) =\spn \left\{ \sum_{\alpha =0}^\beta (-\rmi)^\alpha (\partial_x^\alpha x^\beta)e_\alpha^*
\hspace{2mm}\bigg|\hspace{2mm}  \beta=0, 1,\cdots, m-1  \right \}; \]
\item for $\gamma_{\max} \in I_0:=(-\infty, 1/2)$, the operator \eqref{e:superL} is onto with kernel
\[ \ker (\mathcal{L})= \spn\left\{ \sum_{\alpha =0}^\beta (-\rmi)^\alpha (\partial_x^\alpha x^\beta)e_\alpha 
\hspace{2mm}\bigg|\hspace{2mm} \beta=0, 1,\cdots, m-1 \right\}; \]
\item for $\gamma_{\min} \in I_i$  and $\gamma_{\max} \in I_j$ with $I_k:=(k-1/2, k+1/2)$ for $0< k\in \Z <m$, the kernel of \eqref{e:superL} is
\[ \ker (\mathcal{L})=\spn\left\{ \sum_{\alpha =0}^\beta (-\rmi)^\alpha (\partial_x^\alpha x^\beta)e_\alpha 
\hspace{2mm}\bigg|\hspace{2mm}  \beta=0,1,\cdots,m-j-1 \right\}; \]
and its cokernel is
\[ \cok (\mathcal{L})= \spn \left\{ \sum_{\alpha =0}^\beta (-\rmi)^\alpha (\partial_x^\alpha x^\beta)e_\alpha^* 
\hspace{2mm}\bigg|\hspace{2mm}  \beta=0, 1,\cdots, i-1 \right\}. \]
\end{itemize}
On the other hand, the operator \eqref{e:superL} does not have closed range for $\gamma_-,\gamma_+ \in \{ 1/2, 3/2, \cdots, m-1/2\}$.
\end{Proposition}
The proof of the proposition will occupy the remainder of this section. The key ingredient is the construction of a normal form representation of the operator $L$, through which we conclude that Fredholm properties of the operator $\mathcal{L}$ are equivalent to those of the regularized derivative $[D(-\rmi\partial_x)]^\ell$.  We organize the proof by first establishing Fredholm properties of regularized derivatives defined in the Kondratiev spaces, then Fredholm properties of the normal form of the operator $L$, and eventually concluding the proof by returning to physical space.

\paragraph{Fredholm properties of regularized derivatives.} We employ regularized derivatives as model operators. More specifically, for any $\ell\in\Z^+$ and $\gamma_\pm\in\R$, we define the regularized derivative,
\begin{equation}\label{e:rell2}
\begin{matrix}
[D(-\rmi\partial_x)]^\ell: & \mathcal{D}([D(-\rmi\partial_x)]^\ell)\subset L^2_{\gamma_--\ell, \gamma_+-\ell} & \longrightarrow & L^2_{\gamma_-, \gamma_+} \\
& u & \longmapsto & \partial_x^\ell ( 1 + \partial_x)^{-\ell}u,
\end{matrix}
\end{equation}
with its domain $\mathcal{D}([D(-\rmi\partial_x)]^\ell)=\{u\in L^2_{\gamma_--\ell, \gamma_+-\ell}\mid (1+\partial_x)^{-\ell}u\in M^{\ell,2}_{\gamma_--\ell, \gamma_+-\ell}\}$. Moreover, the Fredholm properties of the operator $[D(-\rmi\partial_x)]^\ell$ are summarized in the following proposition. 
\begin{Proposition}\label{p:regdrv2}
For $\gamma_\pm\in \mathbb{R} \setminus \{ 1/2, 3/2, \cdots, \ell-1/2\}$, the regularized derivative $[D(-\rmi\partial_x)]^\ell$ as defined in \eqref{e:rell} is Fredholm. Moreover, the operator $[D(-\rmi\partial_x)]^\ell$ satisfies the following conditions.
\begin{itemize}
 \item If $\gamma_{max}\in I_0:=(-\infty,1/2)$, the operator $[D(-\rmi\partial_x)]^\ell$ is onto with its kernel equal to $\mathbb{P}_{\ell}(\R)$.
 \item If $\gamma_{min}\in I_\ell:=(\ell-1/2, \infty)$, the operator $[D(-\rmi\partial_x)]^\ell$ is one-to-one with its cokernel equal to $\mathbb{P}_{\ell}(\R)$.
 \item If $\gamma_{min}\in I_i$ and $ \gamma_{max}\in I_j$ with $I_k:=(k-1/2,k+1/2)$ for $0<k\in\Z<\ell$, the kernel and cokernel of the operator $[D(-\rmi\partial_x)]^\ell$ are respectively spanned by $\mathbb{P}_{\ell-j}(\R)$ and $\mathbb{P}_{i}(\R)$.
\end{itemize}
On the other hand, the range of the operator $[D(-\rmi\partial_x)]^\ell$ is not closed if $\gamma_-, \gamma_+\in\{1/2, 3/2,...,\ell-1/2\}$. 
\end{Proposition}
\begin{Proof}
The proof is relegated to Appendix \ref{ss:a1}, where we prove a more general result.
\end{Proof}

\paragraph{Normal form operators.}
We diagonalize every operator $L(k)$ defined in $Y$ into the direct sum of the Fourier counterpart of a regularized derivative and an isomorphism. 
To start with, recalling the definitions of the modified kernel and cokernel expansions  \eqref{e:expmod} and \eqref{e:expmodad}, for any $k\in\Omega_0(C,m)$, we define the projections, 
\begin{equation}
P(k)u:=\langle u, e_0\rangle \tilde{e}(k), \qquad Q(k)v:=\langle v, \tilde{e}^*(k)\rangle e_0^*,
\end{equation}
from which it is straightforward to conclude the following lemma.
\begin{Lemma}\label{l:pro}
There exists $C_0>0$ so that, for any $C>C_0$ and $k\in\Omega_0(C,m)$, the operators
\[
\id - P(k): \langle \tilde{e}(k) \rangle^\perp\rightarrow \langle e_0 \rangle^\perp, \quad 
\id - Q(k): \langle e_0^* \rangle^\perp\rightarrow \langle \tilde{e}^*(k)\rangle^\perp
\] 
are isomorphisms whose inverses take the form,
\begin{equation}
\begin{matrix}
(\id-P(k))^{-1}: & \langle e_0 \rangle^\perp & \longrightarrow &  \langle \tilde{e}(k) \rangle^\perp \\
& u & \longmapsto & u-\frac{\langle u, \tilde{e}(k) \rangle}{\langle \tilde{e}(k), \tilde{e}(k)  \rangle} \tilde{e}(k),
\end{matrix}
\qquad
\begin{matrix}
(\id-Q(k))^{-1}: & \langle \tilde{e}^*(k)\rangle^\perp & \longrightarrow &  \langle e_0^* \rangle^\perp \\
& u & \longmapsto & u-\langle u, e_0^* \rangle e_0^*.
\end{matrix}
\end{equation}
Moreover, for fixed $C>C_0$, both operators and their inverses admit uniform bounds for $k\in\Omega_0(C,m)$.
\end{Lemma}
We also introduce analytic isomorphisms $\iota(k):\langle \tilde{e}(k)\rangle \to \langle e_0^*\rangle$ and $\iota_\perp(k): \langle e_0\rangle^\perp \to \langle \tilde{e}^*(k)\rangle^\perp$. Such isomorphisms can be constructed in many ways and we outline one construction here. Define
\begin{equation}\label{e:pro}
\begin{matrix}
\iota(k):&\langle \tilde{e}(k)\rangle & \longrightarrow & \langle e_0^*\rangle \\
&\alpha \tilde{e}(k) & \longmapsto & \alpha e_0^*,
\end{matrix}
\hspace{2cm}
\begin{matrix}
 \iota_\perp(k): &\langle e_0\rangle^\perp & \longrightarrow & \langle \tilde{e}^*(k)\rangle^\perp\\
&u & \longmapsto & (\id-Q(k))\iota_\perp(0)u,
\end{matrix}
\end{equation}
where we define the isomorphism $\iota_\perp(0): \langle e_0\rangle^\perp \to \langle e_0^*\rangle^\perp$ to be a direct sum of the identity map on 
$\langle e_0\rangle^\perp\cap \langle e_0^* \rangle^\perp$ and a linear length-preserving map from $E_{0,\perp}:=\spn\{e_0^*-\langle e_0^*, e_0 \rangle e_0\}$ 
to $E_{0,\perp}^*:=\spn\{e_0-\langle e_0, e_0^* \rangle e_0^*\}$. More specifically, we choose
\[
\iota_\perp(0)u:=\begin{cases}
u, & u\in \langle e_0\rangle^\perp\cap \langle e_0^* \rangle^\perp,\\
c(e_0-\langle e_0, e_0^* \rangle e_0^*), &u=c(e_0^*-\langle e_0^*, e_0 \rangle e_0)\in E_{0,\perp}.
\end{cases}
\]

We are now ready to define the normal form operators,
\begin{equation}
\begin{matrix}
L_{\rm NF}(k) :&\mathcal{D}(L_{\rm NF}(k))\subset Y &\longrightarrow  &Y \\
 & u & \longmapsto & D^m(k)\iota(k)P(k)u+\iota_\perp(k)(\id - P(k))u,
\end{matrix}
\end{equation}
and prove the following lemma.

\begin{Lemma}[Factorization]\label{l:fact}
For fixed $C>C_0$ and any $k\in\Omega_0(C,m)$, the operator $L(k)$ admits the decomposition,
\[
L(k)=M_\mathrm{L}(k)L_\mathrm{NF}(k)=L_\mathrm{NF}(k)M_\mathrm{R}(k),
\]
where $M_{\rm L\backslash R}:\Omega_0(C, m)\to  B(Y)$ are analytic, $L^\infty$-bounded with an $L^\infty$-bounded inverse.
\end{Lemma}
\begin{Proof}
For $k\neq 0$, the inverse of $L_{\rm NF}(k)$ is analytic and takes the form,
\[
L_{\rm NF}^{-1}(k)u=D^{-m}(k)\iota^{-1}(k)Q(k)u+\iota_\perp^{-1}(k)(\id-Q(k))u=
D^{-m}(k)\langle u, \tilde{e}^*(k) \rangle \tilde{e}(k)+\iota_\perp^{-1}(0)\left( u- \langle u, e_0^*\rangle e_0^*\right).
\]
In addition, we have that, based on \eqref{e:expmod},
\begin{equation*}
\begin{aligned}
\lim_{k\rightarrow 0}L(k)L_{\rm NF}^{-1}(k)u=&
\lim_{k\rightarrow 0}\left[ \frac{(1+\rmi k)^m}{k^m}\langle u, \tilde{e}^*(k)\rangle L(k) \tilde{e}(k)+L(k)\iota^{-1}_\perp(0)\left( u-\langle u, e_0^*\rangle e_0^*\right)\right]\\
=&\lambda_m\langle u, e_0^*\rangle e_0^*+L(0)\iota^{-1}_\perp(0)\left( u-\langle u, e_0^*\rangle e_0^*\right), 
\end{aligned}
\end{equation*}
is an invertible bounded operator.
We now define
\begin{equation}
\label{e:ldef}
M_{\rm L}(k)u:=
\begin{cases}
L(k)L_{\rm NF}^{-1}(k)u, & k\neq 0,\\
\lim_{k\rightarrow 0}L(k)L_{\rm NF}^{-1}(k)u, & k=0,
\end{cases}
\end{equation}
which, according to Riemann's removable singularity theorem and Hypothesis \ref{h:3.2}, implies $M_{\rm L}(k)$ is analytic and invertible for all $k$ in the strip $\Omega_0$. Furthermore, noting that, according to Hypothesis \ref{h:3.1}, $L(k)$ is invertible with uniform bounds for $k\in\Omega_0(C,m)$ with $|\Re k|>k_0$ and
\[
\lim_{\Re k\rightarrow \infty}L_{\rm NF}^{-1}(k)=\langle u, e_0^*\rangle e_0 + \iota_\perp^{-1}(0)\left( u- \langle u, e_0^* \rangle e_0^*\right),
\]
is bounded and invertible, we conclude that $M_{\rm L}(k)$ is uniformly bounded with uniformly bounded inverses. We can define and analyze $M_{\rm R}(k)$ in a completely analogous fashion.
\end{Proof}

\paragraph{Back to physical space --- proof of Proposition \ref{p:f1}.}
We introduce the multiplier operators
\begin{equation}\label{e:multiplier}
\begin{matrix}
\mathcal{M}_{{\rm L}\backslash {\rm R}}: & \mathcal{S}(\R, Y) & \longrightarrow & \mathcal{S}(\R, Y) \\[1mm]
& u(x) & \longmapsto & \widecheck{M_{{\rm L}\backslash {\rm R}}\hat{u}}(x).
\end{matrix}
\end{equation}
which, according to the $L^\infty$-boundedness and invertibility of $\partial_k^\alpha M_{{\rm L}}$ and $\partial_k^\alpha M_{ {\rm R}}$ for all $\alpha\in \Z^+\cup \{0\}$, are isomorphisms on the Schwartz space $\mathcal{S}(\R, Y)$. For any given $\gamma_\pm \in\R$, it is straightforward to see that $\mathcal{S}(\R, Y) \subset L^2_{\gamma_-,\gamma_+}(\R, Y)$ is a continuous embedding. We claim that we can continuously extend the multiplier operators $\mathcal{M}_{{\rm L}\backslash {\rm R}}$ onto $L^2_{\gamma_-,\gamma_+}(\R, Y)$.
In other words, we have the following lemma.
\begin{Lemma}\label{l:L2iso}
For any given $\gamma_\pm \in\R$,  the multiplier operators $\mathcal{M}_{{\rm L}\backslash {\rm R}}: L^2_{\gamma_-,\gamma_+}(\R, Y)\to L^2_{\gamma_-,\gamma_+}(\R, Y)$ are isomorphisms.  
\end{Lemma}

\begin{Remark}
We suspect that results analogous to Lemma \ref{l:L2iso} hold for general anisotropic weighted spaces $L^p_{\gamma_-,\gamma_+}(\R, Y)$ with $p\in(1,\infty)$. It appears however that necessary-and-sufficient condition  for Fourier multipliers on $L^p_{\gamma_-,\gamma_+}(\R, \C)$ with general $p\in(1,\infty)$ are not available, only sufficient conditions such as the Marcinkiewicz and the H\"{o}rmander-Mikhlin multiplier theorems, which both can be generalized to certain families of weighted $L^p(\R,\C)$ spaces; see \cite{muckenhoupt,fefferman1974,kurtz} for details and \cite{amann1, girardiweis, lutz, bu} for general background on operator-valued Fourier multipliers.
\end{Remark}

\begin{Proof}
We first prove the case of isotropic weights, that is, $\gamma_-=\gamma_+=\gamma$. For $\gamma\in\Z_+\cup \{0\}$, we adopt the notation 
$L^2_{\gamma}(\R, Y):=L^2_{\gamma,\gamma}(\R, Y)$ and exploit the Plancherel theorem to derive that
\[
\|\mathcal{M}_{{\rm L}\backslash {\rm R}} u\|_{L^2_{\gamma}(\R, Y)}=\|M_{{\rm L}\backslash {\rm R}}\hat{u}\|_{H^\gamma(\R, Y)}
\leq C(\gamma)\|\hat{u}\|_{H^\gamma(\R, Y)}=C(\gamma)\|u\|_{L^2_\gamma(\R, Y)},
\] 
which, together with a similar inequality for $\mathcal{M}_{{\rm L}\backslash {\rm R}}^{-1}$, shows that $\mathcal{M}_{{\rm L}\backslash {\rm R}}: L^2_{\gamma}(\R, Y)\to L^2_{\gamma}(\R, Y)$ are isomorphisms for $\gamma\in\Z_+\cup \{0\}$ and thus for $\gamma\in \Z_-$ due to duality. By classical interpolation results, see, for example, Theorem 6.4.5 in \cite{bergh}, $H^{n+\theta}(\R, Y)$ is a complex interpolation space between $H^{n}(\R, Y)$ and $H^{n+1}(\R, Y)$ for any given $n\in\Z$ and $\theta\in(0,1)$. Therefore, we conclude that 
$\mathcal{M}_{{\rm L}\backslash {\rm R}}: L^2_{\gamma}(\R, Y)\to L^2_{\gamma}(\R, Y)$ are isomorphisms for $\gamma\in\R$.

To prove the case of anisotropic weights, we start by introducing the exponentially weighted space 
\[
L^2_{\rm exp, \eta}(\R, Y):=\left\{u \in L^1_{\rm loc}(\R, Y)\middle| \rme^{\eta \cdot}u(\cdot)\in L^2(\R, Y)\right\},
\]
with its norm $\|u\|_{L^2_{\rm exp, \eta}(\R, Y)}:=\|\rme^{\eta \cdot}u(\cdot)\|_{L^2(\R, Y)}$ for any given $\eta\in\R$.
Our strategy is to exploit the fact that the space $L^2_{\gamma_-,\gamma_+}(\R, Y)$ admits the decomposition,
\begin{equation}\label{e:intadd}
L^2_{\gamma_-,\gamma_+}(\R, Y)=\left(L^2_{\gamma_-}(\R, Y)\cap L^2_{\rm exp, \eta}(\R, Y)\right) + \left(L^2_{\gamma_+}(\R, Y)\cap L^2_{\rm exp, -\eta}(\R, Y)\right),
\end{equation}
for any $\eta>0$, where norms on intersections and sums are defined in the usual way; see below.
 
With this in mind, we first study the multipliers on $\mathcal{M}_{{\rm L}\backslash {\rm R}}: L^2_{\rm exp, \eta}(\R, Y)\to L^2_{\rm exp, \eta}(\R, Y)$ and claim that they are isomorphisms, for any fixed $|\eta|\leq k_1$, where $k_1$ is half of the width of the strip $\Omega_0(C,m)$. Note that the multiplier on the Schwartz space can be viewed as a convolution operator. More specifically, denoting the reflection $(\mathcal{R}u)(x):=u(-x)$, we define the distribution 
\[
\begin{matrix}
\check{M}_{\rm L\backslash R}: & \mathcal{S}(\R, Y) & \longrightarrow & \C \\
& u & \longmapsto & (\mathcal{M}_{\rm L\backslash R}\mathcal{R}u)(0),
\end{matrix}
\]
from which we readily derive that, for all $u\in\mathcal{S}(\R, Y)$,
\[
(\mathcal{M}_{\rm L\backslash R}u)(x)=(\check{M}_{\rm L\backslash R}\ast u)(x)=\int_\R \check{M}_{\rm L\backslash R}(x-y) u(y)\rmd y.
\]
Since the Fourier transform is given through $\mathcal{F}(\rme^{\eta\cdot}\check{M}_{\rm L\backslash R}(\cdot))(k)=M_{\rm L\backslash R}(k+\rmi \eta)$ for $|\eta|\leq k_1$, we have that the inequality
\[
\begin{aligned}
\|\mathcal{M}_{\rm L\backslash R}u\|_{L^2_{\rm exp, \eta}(\R, Y)}= &\|\int_\R \big[\rme^{\eta(x-y)}\check{M}_{\rm L\backslash R}(x-y) \big]\big[\rme^{\eta y}u(y)\big]\rmd y \|_{L^2(\R, Y)} \\
=& \|\mathcal{F}\big(\rme^{\eta\cdot}\check{M}_{\rm L\backslash R}(\cdot)\big)\mathcal{F}(\rme^{\eta\cdot}u(\cdot)) \|_{L^2(\R, Y)}\\
=& \|M_{\rm L\backslash R}(\cdot+\rmi \eta)\mathcal{F}(\rme^{\eta\cdot}u(\cdot)) \|_{L^2(\R, Y)}\\
\leq&  \|M_{\rm L\backslash R}(\cdot+\rmi \eta)\|_{L^\infty(\R, B(Y))}  \|\mathcal{F}(\rme^{\eta\cdot}u(\cdot)) \|_{L^2(\R, Y)}\\
\leq & C \|u \|_{L^2_{\rm exp, \eta}(\R, Y)},
\end{aligned}
\]
holds for any $|\eta|\leq k_1$ and $u\in \mathcal{S}(\R, Y)$. Noting that $\mathcal{S}(\R, Y)\subset L^2_{\rm exp, \eta}(\R, Y)$ is dense, there are natural extensions of $\mathcal{M}_{\rm L\backslash R}$ as a bounded linear operator on $L^2_{\rm exp, \eta}(\R, Y)$. Analogous reasoning applied to the inverses of $\mathcal{M}_{\rm L\backslash R}$ lets us conclude that the multipliers $\mathcal{M}_{{\rm L}\backslash {\rm R}}: L^2_{\rm exp, \eta}(\R, Y)\to L^2_{\rm exp, \eta}(\R, Y)$ are isomorphisms for any fixed $|\eta|\leq k_1$.

We are now ready to prove the case of anisotropic weights. Given two Banach spaces $E$ and $F$, the linear space $E\cap F$ and $E + F$ are also Banach spaces respectively with norms 
\[
\|u\|_{E\cap F}:= \|u\|_E + \| v \|_F, \quad \| u \|_{E + F} := \inf\{\|v\|_E+\|w\|_F\mid v+w=u, v\in E, w\in F \}.
\] 
Moreover, for a linear operator $L$ bounded on both $E$ and $F$, it is straightforward to check that $L$ is also bounded on $E\cap F$ and $E+F$. Therefore, given $\gamma\pm\in\R$ and $\eta\in[0,k_1]$, due to the fact that $\mathcal{M}_{\rm L\backslash R}$ are isomorphisms on $L^2_{\gamma\pm}$ and $L^2_{\rm exp, \pm\eta}$, we conclude that $\mathcal{M}_{\rm L\backslash R}$ are isomorphisms on the Banach space
\begin{equation}\label{e:bintadd}
B(\gamma_-,\gamma_+, \eta, Y):=\left(L^2_{\gamma_-}(\R, Y)\cap L^2_{\rm exp, \eta}(\R, Y)\right) + \left(L^2_{\gamma_+}(\R, Y)\cap L^2_{\rm exp, -\eta}(\R, Y)\right).
\end{equation}
Es defined in  \eqref{e:intadd}, the Banach spaces $L^2_{\gamma_-,\gamma_+}(\R, Y)$ and $B(\gamma_-,\gamma_+, \eta, Y)$ constitute the same linear space. It is therefore sufficient to show that the natural norm on $L^2_{\gamma_-,\gamma_+}(\R, Y)$ is equivalent to the norm on $B(\gamma_-,\gamma_+, \eta, Y)$ induced by the intersection and sum property. For any $u\in L^2_{\gamma_-,\gamma_+}(\R, Y)$, we have
\[
u=\chi_+u+\chi_-u, \quad \chi_\pm u\in L^2_{\gamma_\pm}(\R, Y)\cap L^2_{\rm exp, \mp\eta}(\R, Y),
\]
and 
\[
\begin{aligned}
\|u\|_{B(\gamma_-,\gamma_+, \eta, Y)}\leq & \| \chi_+u \|_{L^2_{\gamma_+}(\R, Y)\cap L^2_{\rm exp, -\eta}(\R, Y)}+
 \| \chi_-u \|_{L^2_{\gamma_-}(\R, Y)\cap L^2_{\rm exp, \eta}(\R, Y)}\\
 =& \| \chi_+u \|_{L^2_{\gamma_+}(\R, Y)}+\| \chi_+u \|_{L^2_{\rm exp, -\eta}(\R, Y)}+
 \| \chi_-u \|_{L^2_{\gamma_-}(\R, Y)}+\| \chi_-u \|_{ L^2_{\rm exp, \eta}(\R, Y)}\\
 \leq &C(\gamma\pm,\eta)\big[ \| \chi_+u \|_{L^2_{\gamma_+}(\R, Y)}+ \| \chi_-u \|_{L^2_{\gamma_-}(\R, Y)}\big]\\
 =&C(\gamma\pm,\eta)\|u\|_{L^2_{\gamma_-,\gamma_+}(\R, Y)},
 \end{aligned}
\]
which implies that the two norms are equivalent, concluding the proof.
\end{Proof}

Denoting the inverse Fourier transform of $L_{\rm NF}$ as $\mathcal{L}_{\rm NF}$, we have
\[
\mathcal{L}=\mathcal{M}_{\rm L}\mathcal{L}_{\rm NF}, \quad \mathcal{L}^{\rm ad}=\mathcal{M}_{\rm R}^{\rm ad}\mathcal{L}_{\rm NF}^{\rm ad}.
\] The proof of Proposition \ref{p:f1} now reduces to establishing Fredholm properties of $\mathcal{L}_{\rm NF}$. 

\begin{Proof}[of Proposition \ref{p:f1}]
Noting that $Y\cong \langle \tilde{e}(k) \rangle \oplus \langle e_0\rangle^{\perp} \cong\langle e^*_0 \rangle \oplus \langle \tilde{e}^*(k) \rangle^{\perp}$, the normal form operator $L_{\rm NF}(k)$ admits an isomorphic diagonal form, 
\begin{equation}
\begin{matrix}
L_{\rm D}(k):& \langle \tilde{e}(k) \rangle \oplus \langle e_0\rangle^{\perp} &\longrightarrow  &\langle e^*_0 \rangle \oplus \langle \tilde{e}^*(k) \rangle^{\perp} \\
 & \begin{pmatrix} u_1 \\ u_2 \end{pmatrix} & \longmapsto & \begin{pmatrix} D^m(k) \iota(k)& 0 \\ 0 & \iota_{\perp}(k)  \end{pmatrix} \begin{pmatrix} u_1 \\ u_2 \end{pmatrix} 
 \end{matrix}.
\end{equation}
According to Lemma \ref{l:pro}-\ref{l:fact} and definition \eqref{e:pro} of projections $\iota(k)$ and $\iota^\perp(k)$, we derive that
\begin{equation*}
\begin{matrix}
\mathcal{L}_{\rm NF}: & \mathcal{D}(\mathcal{L}_{\rm NF})\subset L^p_{\gamma_--m, \gamma_+-m}(\R, Y) & \longrightarrow & L^p_{\gamma_-, \gamma_+}(R, Y)\\
& u &\longmapsto & \langle D^m(-\rmi \partial_x)u, e_0 \rangle e_0^* +\check{\iota}_\perp ( u-\sum_{j=0}^m\langle D_{C,m}^j(-\rmi \partial_x)u, e_0 \rangle e_j),
\end{matrix}
\end{equation*} 
where $u(x)-\sum_{j=0}^m\langle D_{C,m}^j(-\rmi \partial_x)u(x), e_0 \rangle e_j\in \langle e_0 \rangle^\perp$ for all $x\in\R$ and the mapping
\[
\begin{matrix}
\check{\iota}_\perp: & L^p_{\gamma_-,\gamma_+}(\R, \langle e_0\rangle^\perp) & \longrightarrow & \{u\in L^p_{\gamma_-,\gamma_+}(\R, Y)\mid \sum_{j=0}^m\langle D_{C,m}^j(-\rmi \partial_x)u(x), e_i^*\rangle=0, \text{ for all }x\in\R\}\\[2mm]
& v & \longmapsto &\iota_\perp(0)v-\sum_{j=0}^m\langle D_{C,m}^j(-\rmi \partial_x)[\iota_\perp(0)v], e_i^*\rangle e_0^*,
\end{matrix}
\]
is an isomorphism. 
As a result, Fredholm properties of $\mathcal{L}_{\rm NF}$ are encoded in the regularized derivative operator $[D(-\rmi \partial_x)]^m$. More specifically, we note that
\[
\mathcal{F}^{-1}\left[ D^m(k)\iota(k)\Big(\hat{u}(k) \tilde{e}(k)\Big) \right]
=\Big( [D(-\rmi \partial_x)]^m u(x)\Big) e_0^*, \quad
\mathcal{F}^{-1}\Big(\hat{u}(k) \tilde{e}(k)\Big)=\sum_{j =0}^m \Big(\left[D_{C,m} (-\rmi \partial_x)\right]^j u(x) \Big)e_j,
\]
which implies that the kernel and cokernel of $\mathcal{L}_{\rm NF}$ are given respectively by
\begin{equation*}
\begin{aligned}
&\ker(\mathcal{L}_{\rm NF}) = \left \{ \sum_{j =0}^m \Big(\left[D_{C,m} (-\rmi \partial_x)\right]^j u(x) \Big)e_j
 \hspace{2mm}\bigg|\hspace{2mm} u(x) \in \ker \Big( [D(-\rmi \partial_x)]^m \Big)  \right\},\\
 &\cok(\mathcal{L}_{\rm NF}) = \left \{ \sum_{j =0}^m \Big(\left[\overline{D_{C,m} (\rmi \partial_x)}\right]^j u(x) \Big)e_j^*
 \hspace{2mm}\bigg|\hspace{2mm} u(x) \in \cok \Big( [D(-\rmi \partial_x)]^m \Big)  \right\} .
 \end{aligned}
\end{equation*}
The statements in Proposition \ref{p:f1} then follow by applying the statement of Proposition \ref{p:regdrv} to the above analysis and noting that, for any $u\in\mathbb{P}_m(\R)$,
\[
\left[D_{C,m} (-\rmi \partial_x)\right]^j u(x)=(-\rmi)^\alpha\partial_x^\alpha u(x).
\]
\end{Proof}

\subsection{Operators with discrete translation symmetry}\label{s:3.2}

The results from Section \ref{s:3.1} can be easily adapted to the case of an operator $\mathcal{L}$ on $\ell^2(\Z, Y)$, that commutes with the discrete translation group $\Z$. The discrete Fourier transform takes the form
\begin{equation}
\label{e:disft}
\begin{matrix}
\mathcal{F}_d: & \ell^2(\Z, Y) & \longrightarrow &L^2(\mathcal{T}_1, Y)\\
&\underline{u}= \{u_j\}_{j\in\Z}& \longmapsto & \hat{u}(\sigma)=\sum_{j\in\Z}u_j\rme^{-2\pi\rmi j\sigma},
\end{matrix}
\end{equation}
where $\mathcal{T}_1:=\R/\Z$ denotes the unit circle.
The counterparts of the derivative $\partial_x$ are the discrete derivatives,
\begin{equation}
\delta_+ (\{a_j\}_{j\in\Z}) := \{a_{j+1}-a_{j}\}_{j\in\Z},\hspace{0.6cm}
\delta_- (\{a_j\}_{j\in\Z}) := \{a_j-a_{j-1}\}_{j\in\Z}, \hspace{0.6cm}
\delta:=-\rmi(\delta_++\delta_-)/2.
\end{equation}
The Fourier transform of $\mathcal{L}$, denoted as $\hat{\mathcal{L}}=\int_{\mathcal{T}_1}L(\sigma)\rmd \sigma$, is an isomorphism of $L^2(\mathcal{T}_1,Y)$, that is,
\begin{equation}
\begin{matrix}
\hat{\mathcal{L}}: & \mathcal{D}(\hat{\mathcal{L}})\subset L^2(\mathcal{T}_1,Y) & \longrightarrow & L^2(\mathcal{T}_1,Y) \\
& u(\sigma) & \longmapsto & L(\sigma)u(\sigma),
\end{matrix}
\end{equation}
with $L(\sigma)$ linear and bounded on $Y$ for all $\sigma\in\mathcal{T}_1$. 
\begin{Hypothesis}[Analyticity, periodicity and simple kernel]\label{h:3.3}
We assume that $L(\sigma)$ is analytic, uniformly bounded, $1$-periodic, with values in the set of bounded operators on $Y$,  in a strip $\sigma\in \Omega_1:=\R\times (-\rmi \sigma_\rmi,\rmi \sigma_\rmi)$ for some $\sigma_\rmi>0$. 
Moreover, we require that $L(\sigma)$, restricted to $\sigma\in[-1/2, 1/2]$, is invertible except at $\sigma=0$ and $L(0)$ admits a simple kernel spanned by $e_0$ with $\langle e_0, e_0 \rangle=1$.
\end{Hypothesis}

\begin{Remark}
For convenience, we identify the interval $[-1/2, 1/2]$ with the unit circle $\mathcal{T}_1$, collapsing endpoints $-1/2\sim 1/2$.
\end{Remark}

We adopt all the notations in the continuous case, except for those related to pseudo-derivative symbols.
The new pseudo-derivatives take the following forms,
\begin{equation}\label{e:dispd}
D_+(\sigma)=\rme^{2\pi\rmi\sigma}-1, \qquad D_-(\sigma)=1-\rme^{-2\pi\rmi\sigma}, \qquad 
D_{C,m}(\sigma)=(\rme^{2\pi\rmi\sigma}-1)\left[1+\rmi C\sin^m(2\pi\sigma)\right]^{-1},
\end{equation}
whose associated physical operator are respectively $\delta_+$, $\delta_-$ and $\delta_+\left[1+\rmi C\delta^m\right]^{-1}$.
Here $m\in\Z^+$ is the power related to the expansion of the zero eigenvalue, $\lambda(\sigma)=\lambda_m \sigma^m+O(\sigma^{m+1})$, with $\lambda_m\neq0$ for $\sigma\sim 0\in\C$. 
The constant $C>0$ will eventually be chosen sufficiently large so that the norm of the bounded multiplier $D_{C,m}$ is arbitrarily small. 
As a matter of fact, in the strip
\[
\Omega_1(C,m):=\left\{\sigma\in\Omega_1\middle| {|\Re \sigma|\leq 1/2,} |\Im \sigma|< \frac{1}{2\pi}\sinh^{-1}\left(\frac{1}{\sqrt[m]{2C}}\sin(\frac{\pi}{2m})\right)\right\},
\] 
$D_{C,m}(\sigma)$ is analytic and uniformly bounded, that is, there exists a constant C(m) so that
\[
\|D_{C,m}(\sigma)\|\leq \frac{C(m)}{\sqrt[m]{C}}, \text{ for all }\sigma\in\Omega_1(C,m).
\]

Moreover, we define $e(\sigma)=\sum_{j=0}^{m} e_j\sigma^j$ and $e^*(\sigma)=\sum_{j=0}^{m} e_j^*\bar{\sigma}^j$ so that 
\begin{equation*}
L(\sigma)e(\sigma)=\rmO(\sigma^m),\quad L^*(\sigma)e^*(\sigma)=\rmO(\sigma^m), \quad 
\left\langle \sum_{j=0}^{m-1} L_{m-j}e_j,e_0^*\right\rangle\neq0, \quad
\sum_{j=0}^k L_je_{k-j}=0,\quad k=0,\ldots,m-1.
\end{equation*}
 There exist $\{\tilde{e}_j, \tilde{e}_j^*\}_{j=0}^m\subset Y$, independent of $C$, and 
\[
\tilde{e}(\sigma):=\sum_{j=0}^{m} \left[D_{C,m}(\sigma)\right]^j\tilde{e}_j, \qquad
\tilde{e}^*(\sigma):=\sum_{j=0}^{m} \left[\overline{D_{C,m}(\sigma)}\right]^j\tilde{e}_j^*, \qquad
\sigma\in\Omega_1(C,m).
\]
so that $L(\sigma)\tilde{e}(\sigma)=\rmO(\sigma^{m})$ and $L^*(\sigma)\tilde{e}^*(\sigma)=\rmO(\sigma^{m})$.

\begin{Proposition}[Fredholm properties of $\mathcal{L}$]\label{p:f2}
For  $\gamma_\pm \not\in \{1/2,3/2, \cdots, m-1/2\}$, the operator satisfying Hypothesis \ref{h:3.3},
\begin{equation}\label{e:dissuperL}
\mathcal{L}:\mathcal{D}(\mathcal{L})\subset \ell^2_{\gamma_--m,\gamma_+-m}(\Z,Y)\to  \ell^2_{\gamma_-,\gamma_+}(\Z,Y),
\end{equation}
 is closed, densely defined, and Fredholm.
Letting $\gamma_{\max} = \max \{ \gamma_-, \gamma_+\} , \gamma_{\min} = \min\{ \gamma_-, \gamma_+\}$ and $\underline{\eta}^\beta:=\{\eta^\beta\}_{\eta\in\Z}$, we have that 
\begin{itemize}
\item for $\gamma_{\min} \in I_m:=(m-1/2, \infty)$, the operator \eqref{e:dissuperL} is one-to-one with cokernel
\[ \cok =\spn \left\{ \sum_{\alpha =0}^\beta (\delta_+^\alpha \underline{\eta}^\beta)\tilde{e}_\alpha^* 
 \hspace{2mm}\bigg|\hspace{2mm} \beta=0, 1,\cdots, m-1 \right\}; \]
\item for $\gamma_{\max} \in I_0:=(-\infty, 1/2)$, the operator \eqref{e:dissuperL} is onto with kernel
\[ \ker =\spn \left\{ \sum_{\alpha =0}^\beta (\delta_+^\alpha \underline{\eta}^\beta)\tilde{e}_\alpha 
 \hspace{2mm}\bigg|\hspace{2mm} \beta=0, 1,\cdots, m-1 \right\}; \]
\item for $\gamma_{\min} \in I_i$ and $\gamma_{\max} \in I_j$ with $I_k:=(k-1/2, k+1/2)$ for $0<k\in\Z<m$, the kernel of \eqref{e:dissuperL} is
\[ \ker =  \spn \left\{ \sum_{\alpha =0}^\beta (\delta_+^\alpha \underline{\eta}^\beta)\tilde{e}_\alpha 
\hspace{2mm}\bigg|\hspace{2mm} \beta=0, 1,\cdots, m-j-1 \right\};  \]
and its cokernel is
\[ \cok =\spn \left\{ \sum_{\alpha =0}^\beta (\delta_+^\alpha \underline{\eta}^\beta)\tilde{e}_\alpha^*  
 \hspace{2mm}\bigg|\hspace{2mm} \beta=0, 1,\cdots, i-1 \right\}.\]
\end{itemize}
On the other hand, the operator \eqref{e:dissuperL} does not have closed range for $\gamma_-,\gamma_+ \in \{ 1/2, 3/2, \cdots, m-1/2\}$.
\end{Proposition}
\begin{Proof}
Just as in  the continuous case, the proof reduces to the verification of Fredholm properties of the discrete derivative $\delta_+^{m-j}\delta_-^{j}$, for $j=0,1,\cdots, m$, which is relegated to Appendix \ref{ss:a2}.
\end{Proof}

\subsection{Floquet-Bloch theory and periodic coefficients}\label{s:3.3}

We are interested in operators posed on the real line, with only a discrete translational symmetry. Examples are of course the linearization at periodic structures, but include more generally operators with periodic coefficients, $\mathcal{P}(\partial_x,x)$, periodic in $x$. One commonly introduces 
the Bloch-wave transform
\begin{equation*}
\label{e:bloch}
 \begin{matrix}
  \mathcal{B}: &L^2(\mathcal{T}_1, [L^2([0,2\pi])]^n)&\longrightarrow &[L^2(\R)]^n\\
  &\mathbf{U}(\sigma,x)&\longmapsto& \int_{\mathcal{T}_1}\rme^{\rmi\sigma x}\mathbf{U}(\sigma,\cdot)\rmd\sigma,
 \end{matrix}
\end{equation*}
which is an isometric isomorphism with its inverse
\begin{equation}
\label{e:bin}
\begin{matrix}
 \mathcal{B}^{-1}: &[L^2(\R)]^n &\longrightarrow&L^2(\mathcal{T}_1, [L^2([0,2\pi])]^n)\\
 &\mathbf{u}(x)&\longmapsto&\frac{1}{2\pi}\sum_{\ell\in \Z}\rme^{\rmi\ell x}\widehat{\mathbf{u}}(\sigma+\ell).
\end{matrix}
\end{equation} 
We refer to \cite[\Rmnum{13}.16.]{reedsimon} for details. Under the Bloch-wave transform, $\mathcal{P}(\partial_x,x)$ defined on $[L^2(\R)]^n$ becomes a direct integral --- the Bloch-wave decomposition,
\begin{equation}
\label{e:blod}
\mathcal{B}^{-1}\circ \mathcal{P} \circ\mathcal{B}=\int_{\mathcal{T}_1} P_{\mathrm{BL}}(\sigma)\rmd \sigma,
\end{equation}
where the Bloch-wave operator $P_{\mathrm{BL}}(\sigma)$ takes the form
\begin{equation}
\label{e:blow}
\begin{matrix}
P_{\mathrm{BL}}(\sigma): & \mathcal{D}(P_{\mathrm{BL}}(\sigma))\subset [L^2([0,2\pi])]^n & \longrightarrow & [L^2([0,2\pi])]^n\\
& u(x) & \longmapsto & P(\partial_x+\rmi\sigma, x)u(x).
\end{matrix}
\end{equation}
We assume that the family of Bloch-wave operators $P_{\rm BL}(\sigma)$ satisfies the following hypothesis.
\begin{Hypothesis}[Analyticity and simple kernel]\label{h:3.4}
We assume that $P_{\rm BL}(\sigma)$ is analytic and  uniformly bounded, 1-periodic, with values in the set of bounded operators on $Y$,  in a strip $\sigma\in \Omega_1:={\R}\times (-\rmi \sigma_\rmi,\rmi \sigma_\rmi)$ for some $\sigma_\rmi>0$. 
Moreover, we require that $P_{\rm BL}(\sigma)$,restricted to $[-1/2,1/2]$, is invertible except at $\sigma=0$ and $P_{\rm BL}(0)$ admits a simple kernel spanned by $e_0$ with $\langle e_0, e_0 \rangle=1$.
\end{Hypothesis}

In order to exploit the results from Section \ref{s:3.2}, we first define the chopping operator $\mathcal{C}$ that identifies $[L^2(\R)]^n$ with $\ell^2(\Z, [L^2([0,2\pi])]^n)$, that is,
\begin{equation*}
\label{e:chop}
\begin{matrix}
\mathcal{C}: & [L^2(\R)]^n & \longrightarrow & \ell^2(\Z, [L^2([0,2\pi])]^n)\\
& u& \longmapsto & \{u(2\pi j+x)\}_{j\in\Z},
\end{matrix}
\end{equation*}
and the discrete Fourier transform taking the form
\begin{equation}
\label{e:fourier}
\begin{matrix}
\mathcal{F}_d: & \ell^2(\Z, [L^2([0,2\pi])]^n) & \longrightarrow &L^2(\mathcal{T}_1, [L^2([0,2\pi])]^n)\\
&\underline{u}= \{u_j\}_{j\in\Z}& \longmapsto & \sum_{j\in\Z}u_j(x)\rme^{-2\pi\rmi j\sigma}.
\end{matrix}
\end{equation}
Under the transformations $\mathcal{C}$ and $\mathcal{F}_{\rm d}$, $\mathcal{P}(\partial_x,x)$ again becomes a direct integral with the notation
\begin{equation}
\label{e:abloch}
\int_{\mathcal{T}_1} P(\sigma)\rmd \sigma:=\mathcal{F}_{\rm d}\circ\mathcal{C}\circ \mathcal{P} \circ \mathcal{C}^{-1}\circ \mathcal{F}_{\rm d}^{-1}.
\end{equation}
In fact, for any $U\in \mathcal{D}(\int_{\mathcal{T}_1} P(\sigma)\rmd \sigma)$, we have that 
\begin{equation*}
\begin{aligned}
\Big(\mathcal{F}_{\rm d}\circ\mathcal{C}\circ \mathcal{P} \circ \mathcal{C}^{-1}\circ \mathcal{F}_{\rm d}^{-1}(U)\Big)(\sigma,x)
&=\sum_{j\in\Z}\rme^{-2\pi\rmi j\sigma}\left(\mathcal{P}(\partial_x,x)\int_{\mathcal{T}_1} U(\eta, x)\rme^{2\pi\rmi j \eta}\rmd \eta\right)\\
&=\mathcal{P}(\partial_x,x)\int_{\mathcal{T}_1} U(\eta,x)\left(\sum_{j\in\Z}\rme^{2\pi\rmi j (\eta-\sigma)}\right)\rmd \eta\\
&=\mathcal{P}(\partial_x,x)\int_{\mathcal{T}_1} U(\eta,x)\delta(\eta-\sigma)\rmd \eta\\
&=\mathcal{P}(\partial_x,x)U(\sigma,x),
\end{aligned}
\end{equation*}
which shows that, for any $\sigma\in \mathcal{T}_1$,
\begin{equation*}
\begin{matrix}
P(\sigma): & \mathcal{D}(P(\sigma))\subset [L^2([0,2\pi])]^n & \longrightarrow & [L^2([0,2\pi])]^n\\
& u(x) & \longmapsto & \mathcal{P}(\partial_x, x)u(x).
\end{matrix}
\end{equation*}
We conclude with a commutative diagram of isomorphisms as follows, dropping the superscript $n$ for ease of notation,
\begin{equation*}
\label{e:cmd}
\begin{matrix}
L^2(\mathcal{T}_1, L^2([0,2\pi]))&\overset{\mathcal{B}}{\longrightarrow}&
L^2(\R)&\overset{\mathcal{C}}{\longrightarrow}&
\ell^2(\Z, L^2([0,2\pi]))&\overset{\mathcal{F}_{\rm d}}{\longrightarrow}&
L^2(\mathcal{T}_1, L^2([0,2\pi]))\\
\downarrow\int_{\mathcal{T}_1} P_{\mathrm{BL}}(\sigma)\rmd \sigma&&\downarrow \mathcal{P} &&&&\downarrow \int_{\mathcal{T}_1} P(\sigma)\rmd \sigma\\
L^2(\mathcal{T}_1, L^2([0,2\pi]))&\overset{\mathcal{B}}{\longrightarrow}&
L^2(\R)&\overset{\mathcal{C}}{\longrightarrow}&
\ell^2(\Z, L^2([0,2\pi]))&\overset{\mathcal{F}_{\rm d}}{\longrightarrow}&
L^2(\mathcal{T}_1, [^2([0,2\pi])),
\end{matrix}
\end{equation*}
from which it is straightforward to see that $\int_{\mathcal{T}_1} P_{\mathrm{BL}}(\sigma)\rmd \sigma$ and $\int_{\mathcal{T}_1} P(\sigma)\rmd \sigma$ are isomorphic. Moreover, we have the following lemma.
\begin{Lemma}\label{l:bceq}
The operators $P(\sigma)$ and $P_\mathrm{BL}(\sigma)$ are canonically isomorphic for all $\sigma\in\mathcal{T}_1$.
\end{Lemma}
\begin{Proof}
From (\ref{e:bin}-\ref{e:blod}) and (\ref{e:fourier}-\ref{e:abloch}), we summarize that for any $\sigma\in \mathcal{T}_1$, 
\[
\mathcal{D}(P(\sigma))=\{\rme^{\rmi\sigma x}u(x)\in [L^2([0,2\pi])]^n\mid  u(x)\in\mathcal{D}(P_{\mathrm{BL}}(\sigma))\},
\]
which directly implies that we have the isomorphism
\begin{equation}\label{e:isoblo}
P_{\mathrm{BL}}(\sigma)=\rme^{-\rmi \sigma x}P(\sigma)\rme^{\rmi \sigma x}.
\end{equation}
\end{Proof}

According to Hypothesis \ref{h:3.4}, there exist $m\in\Z^+$, $\lambda_m\neq0$, $e(\sigma)=\sum_{j=0}^{m} e_j\sigma^j$ and $e^*(\sigma)=\sum_{j=0}^{m} e_j^*\bar{\sigma}^j$ with 
\begin{equation}\label{e:expmod2}
P_{\rm BL}(\sigma)e(\sigma)=
\lambda_me_0\sigma^m+\rmO(\sigma^{m+1}), 
\end{equation}
and
\begin{equation}\label{e:expmodad2}
P_{\rm BL}^*(\sigma)e^*(\sigma)=
\bar{\lambda}_me_0^*\sigma^m+\rmO(\sigma^{m+1}),
\end{equation}
so that 
\begin{equation*}
\left\langle \sum_{j=0}^{m-1} P_{{\rm BL}, m-j}e_j,e_0^*\right\rangle\neq0, \quad
\sum_{j=0}^k P_{{\rm BL}, j}e_{k-j}=0,\quad k=0,\ldots,m-1.
\end{equation*}
According to Lemma \ref{l:bceq} and Proposition \ref{p:f2}, we have the following proposition. 
\begin{Proposition}[Fredholm properties of $\mathcal{L}$]\label{p:f3}
For  $\gamma_-,\gamma_+ \not\in \{1/2,3/2, \cdots, m-1/2\}$, the operator satisfying Hypothesis \ref{h:3.4},
\begin{equation}\label{e:dissuperL2}
\mathcal{P}:\mathcal{D}(\mathcal{P})\subset L^2_{\gamma_--m,\gamma_+-m}\to  L^2_{\gamma_-,\gamma_+},
\end{equation}
 is closed, densely defined, and Fredholm.
 Letting $\gamma_{\max} = \max \{ \gamma_-, \gamma_+\} , \gamma_{\min} = \min\{ \gamma_-, \gamma_+\}$, we have that 
\begin{itemize}
\item for $\gamma_{\min} \in I_m:=(m-1/2, \infty)$, the operator \eqref{e:dissuperL2} is one-to-one with cokernel
\[ \cok =\spn \left\{  \sum_{\alpha =0}^\beta \frac{(\rmi x)^\alpha}{\alpha !}e_{\beta-\alpha}^*
 \hspace{2mm}\bigg|\hspace{2mm} \beta=0, 1,\cdots, m-1 \right\}; \]
\item for $\gamma_{\max} \in I_0:=(-\infty, 1/2)$, the operator \eqref{e:dissuperL2} is onto with kernel
\[ \ker =\spn \left\{ \sum_{\alpha =0}^\beta \frac{(\rmi x)^\alpha}{\alpha !}e_{\beta-\alpha}
 \hspace{2mm}\bigg|\hspace{2mm} \beta=0, 1,\cdots, m-1 \right\}; \]
\item for $\gamma_{\min} \in I_i$ and $\gamma_{\max} \in I_j$ with $I_k:=(k-1/2, k+1/2)$ for $0<k\in\Z<m$, the kernel of \eqref{e:dissuperL2} is
\[ \ker =  \spn \left\{  \sum_{\alpha =0}^\beta \frac{(\rmi x)^\alpha}{\alpha !}e_{\beta-\alpha}
 \hspace{2mm}\bigg|\hspace{2mm} \beta=0, 1,\cdots, m-j-1 \right\};  \]
and its cokernel is
\[ \cok =\spn \left\{  \sum_{\alpha =0}^\beta \frac{(\rmi x)^\alpha}{\alpha !}e_{\beta-\alpha}^*
 \hspace{2mm}\bigg|\hspace{2mm} \beta=0, 1,\cdots, i-1 \right\}.\]
\end{itemize}
On the other hand, the operator \eqref{e:dissuperL2} does not have closed range for $\gamma_-,\gamma_+ \in \{ 1/2, 3/2, \cdots, m-1/2\}$.
\end{Proposition}

\begin{Proof}
All results in this proposition, except explicit forms of kernels and cokernels, are direct consequences of Proposition \ref{p:f2}. From the isomorphism property  \eqref{e:isoblo} and the expansion \eqref{e:expmod2}, we have, for $\beta=0, 1, \cdots, m-1$,
\[
\mathcal{P} \sum_{\alpha =0}^\beta \frac{(\rmi x)^\alpha}{\alpha !}e_{\beta-\alpha} =0,
\]
which, combined with the domain of $\mathcal{P}$ for given $\gamma_\pm$, concludes the proof.
\end{Proof}

\begin{Remark}
There is an alternative  way to obtain the explicit forms of kernels and cokernels. The first step is to obtain explicit forms of $\tilde{e}_j$ and $\tilde{e}_j^*$. Taking $\tilde{e}_j$ for example,  we note that the first $m+1$ terms of the Taylor expansion of $\rme^{\rmi x\sigma}e(\sigma)$ and $\sum_{j=0}^m(\rme^{2\pi\rmi\sigma}-1)^j \tilde{e}_j$ with respect to $\sigma$ are the same. More specifically, we have 
\begin{equation*}
\begin{aligned}
\rme^{\rmi x\sigma}e(\sigma)&=e_0+\sum_{k=1}^m\left(\sum_{j=0}^k\frac{(\rmi x)^j}{j!}e_{k-j}\right)\sigma^k+O(\sigma^{m+1}),\\
\sum_{j=0}^m(\rme^{2\pi\rmi\sigma}-1)^j \tilde{e}_j&=\tilde{e}_0+\sum_{k=1}^m\frac{(2\pi\rmi)^k}{k!}\left(A(k,j)\tilde{e}_j\right)\sigma^k+O(\sigma^{m+1}),
\end{aligned}
\end{equation*}
where 
\[
A(k,j)=\sum_{\ell=1}^j\binom{j}{\ell}\ell^k(-1)^{j-\ell},
\]
with $A(k,j)=0$ for $1<k<j$. We can then solve $\{\tilde{e}_j\}_{j=0}^{m}$ in terms of $\{e_j\}_{j=0}^{m}$. In a second step, we plug all these explicit expansions of $\tilde{e}_j$'s into Proposition \ref{p:f2} to derive explicit forms of kernels and cokernels.
\end{Remark}




\section{Impurities}\label{s:4}

We now prove Theorem \ref{t:1}. Recalling $\chi_\pm$ is a smooth partition of unity with $\mathrm{supp}(\chi_+)\subset (-1,\infty)$, $\chi_-(x)=\chi_+(-x)$, we write $\theta=\chi_+-\chi_-$ and
\begin{equation}\label{e:phi}
\begin{aligned}
\varphi(x)&=k_0 x-\varphi_0 + k_1 \Theta- \varphi_1 \theta(x),& \varphi'(x)&=k_0 + k_1\theta(x) -\varphi_1\theta'(x),\\
\varphi^\pm(x)&=k_0 x-\varphi_0 \pm(k_1 x-\varphi_1), & (\varphi^\pm)'(x)&=k_0\pm k_1,
\end{aligned}
\end{equation}
where $\Theta(x):=\int_0^x\theta(y)\rmd y+c$ with the constant $c>0$ chosen so that $\Theta(x)=|x|$ for $|x|>1$.
We think of $\varphi_j$ and $k_j$ as matching variables in the far field and we will consider $\psi_0=(\varphi_0,k_0)$ as free parameters  and $\psi_1=(\varphi_1,k_1)$ as variables, and write $\psi=(\psi_0,\psi_1)$, so that $\varphi=\varphi(x;\psi),\varphi^\pm=\varphi^\pm(x;\psi) $. We write 
\begin{equation}\label{e:defup}
\up^\psi(x):=\up(k_*x+\varphi(x;\psi);k_*+\varphi'(x;\psi)),\qquad \up^{\pm,\psi}(x):=\up(k_*+\varphi^\pm(x;\psi);k_*+(\varphi^\pm)'(x;\psi)).
\end{equation}

We then substitute the ansatz $u(x)=\up^\psi+w$ into the stationary Swift-Hohenberg equation, to obtain
\begin{equation}\label{e:ansatz}
L_\mathrm{SH}(\up^\psi+w)+F(\up^\psi+w)+\varepsilon g=0,
\end{equation}
where 
\[
L_\mathrm{SH}=-(1+\partial_x^2)^2,\qquad F(u)=\mu u - u^3.
\]
The phase shifts $\varphi^\pm$ encode simply shifted phases and wavenumbers, so that $\up^{\pm,\psi}$ are solutions to the Swift-Hohenberg equation and, for both $+$ and $-$, 
\[
\chi_\pm\left(L_\mathrm{SH} \up^{\pm,\psi}+F(\up^{\pm,\psi})\right)=0.
\]
Subtracting these from \eqref{e:ansatz} gives
\begin{equation}\label{e:bif}
L_\mathrm{SH} w+F'(\up^\psi)w+N(w, \psi)+K+\varepsilon G=0,
\end{equation}
where
\[
N(w,\psi)=F(\up^\psi+w)-F(\up^\psi)-F'(\up^\psi)w=\rmO(w^2),\quad G=g(x,\up^\psi+w),
\]
and the commutator $K$ depends on $\psi$, only,
\[
K=L_\mathrm{SH}\up^\psi-\sum_\pm\chi_\pm L_\mathrm{SH} \up^{\pm,\psi}+F(\up^\psi)-\sum\chi_\pm F(\up^{\pm,\psi}).
\]
In particular, one readily finds that $K$ is compactly supported and smooth in $\psi$ as an element of $H^k_\gamma$ for any $k,\gamma$. Expanding 
\[
K=K_1\cdot \psi+ K_2, \quad K_2= \rmO(|\psi|^2),
\]
gives 
\begin{equation}\label{e:bif1}
\mathcal{L}^\psi (w,\psi)+\mathcal{N}(w,\psi)+\varepsilon G(w,\psi)=0,
\end{equation}
where 
\[
\mathcal{L}^\psi (w,\psi)=L_\mathrm{SH} w+F'(\up^\psi)w+K_1\cdot\psi,
\]
 with the following notation
 $$
 K_1:=\partial_\psi K|_{\psi=0}=(K_{\varphi_0},K_{k_0},K_{\varphi_1},K_{k_1}), \quad \mathcal{N}(w,\psi):=N(w,\psi)+K_2=\rmO(|w|^2+|\psi|^2).
 $$
Our goal is to use Lyapunov-Schmidt reduction to solve \eqref{e:bif1} with variables $w,\psi_1$ and parameters $\varepsilon,\psi_0$, near the trivial solution $k_0=k_1=\varphi_1=\varepsilon=0$, $w=0$,  and fixed $\varphi_0\in [0,2\pi)$. 
\begin{Remark}\label{r:shift}
Without loss of generality, we can also redefine the primary pattern, shifting its location by $\frac{\varphi_0}{k_*}$ in a $\varphi_0$-dependent fashion, and subsequently applying the shift $x'=x-\frac{\varphi_0}{k_*}$ in \eqref{e:sh}. As a consequence, in our proof, $\varphi_0\equiv 0$, or, in other words, $\varphi_0$ as a variable does not appear within $\up^\psi$ and the dependence on $\varphi_0$ is moved to $g=g(x'+\frac{\varphi_0}{k_*},u)$.
\end{Remark}
Making the role of variables versus parameters explicit, we further decompose 
\[
\mathcal{L}^\psi (w,\psi)=\mathcal{L}_1^\psi (w,\psi_1)+\mathcal{L}_0^\psi \psi_0,
\]
with
\[ \mathcal{L}^\psi_1(w,\psi_1) = L_\mathrm{SH}w+ F'(\up^\psi)w + K_{\varphi_1} \varphi_1 + K_{k_1}k_1, \qquad \mathcal{L}_0^{\psi}\psi_0 =  K_{\varphi_0} \varphi_0+ K_{k_0}k_0. \]

In order to implement Lyapunov-Schmidt reduction, we proceed as follows. We precondition \eqref{e:bif1} with $\mathcal{M}(\psi):=(\mathcal{L}_1^\psi)^{-1}$ and consider the resulting equation 
\[
(w,\psi_1)+\mathcal{M}(\psi)\left(\mathcal{L}_0^\psi \psi_0+\mathcal{N}(w,\psi)+\varepsilon G(w,\psi)\right)=0,
\]
on $H^4_{\gamma_*-3-\delta}\times \R^2$, in a neighborhood of the origin, with parameters $\psi_0,\varepsilon$. The following two ingredients ensure that we can actually apply the implicit function theorem near the trivial solution $w=\psi_1=0$. 
\begin{enumerate}
\item The inverse 
$\mathcal{M}(\psi)$ is bounded from $L^2_\gamma$ to $H^4_{\gamma-2}\times \R^2$, and $C^1$ in $\psi$ when considered as an operator from $L^2_\gamma$ to $H^4_{\gamma-3-\delta}$, for $\gamma>3/2$.
\item The nonlinearity $\mathcal{N}$ is of class $C^1$ as a map from $H^4_{\gamma}\times \R^4$ into $L^2_{2\gamma}$, with vanishing derivatives at the origin. 
\end{enumerate}
We then choose $\gamma=\gamma_*$ in (i) and $2\gamma=\gamma_*$ in (ii), which gives the restriction $2(\gamma_*-3-\delta)>\gamma_*$, compatible with $\gamma_*>6$. 

The second part is quite standard, using that $u\mapsto u\cdot u$ maps $H^k_\gamma$ into $H^k_{2\gamma}$ for $k>1/2$, and we will focus on the first part in the next two sections. 
We therefore proceed in several steps. We first show bounded invertibility for $\psi=0$ in section \ref{s:4.1} , in particular computing the derivatives of $K$ and their projection on the cokernel of $ \mathcal{L}^0_1= L_\textrm{SH} + F'(\up)$, where $\up$ simply stands for $\up(\xi; k_*)$. We then show bounded invertibility and continuity of $\mathcal{L}^{\psi}_1$ for  $\psi\neq 0$ using a decomposition argument in Section \ref{s:4.2}. Finally, we compute expansions in Section \ref{s:4.3}. 

\subsection{Invertibility at $\psi \equiv 0$}\label{s:4.1}

In this subsection we drop the subscripts from $\mathcal{L}^0_1$. We first show that 
\begin{equation}\label{e:L0}
 \mathcal{L}^0=L_\textrm{SH} + F'(\up),
\end{equation} 
is Fredholm and identify the cokernel, then compute projections of the partial derivatives of $K_1$ on the cokernel, and finally identify projection coefficients with effective diffusivity.
Recall that  $\up(\xi;k_*)$, with $\xi = k_* x$, denotes a periodic solution to the unperturbed Swift-Hohenberg equation. Throughout this section we will write $\up':=\partial_x \up=k_*\partial_\xi \up(\xi; k_*)$, $\partial_\xi \up:=\partial_\xi \up(\xi; k_*)$ and $\partial_k\up:=\partial_k\up(\xi; k_*)$.

\paragraph{Fredholm properties of $\mathcal{L}^0$.}

We start by putting the results from Section \ref{s:3} to work.

\begin{Proposition}\label{prop:fredholm}
Assume Hypotheses \ref{h:0}--\ref{h:2}. For all $\gamma>3/2$, the linear operator $\mathcal{L}^0:\mathcal{D}(\mathcal{L}^0)\subset H^4_{\gamma-2}\to L^2_\gamma$ is Fredholm of index -2, with trivial kernel and cokernel spanned by $\up'$ and $\upk=x\partial_\xi\up+\partial_k\up$. 
\end{Proposition}
\begin{Proof}
According to Proposition \ref{p:f3} and the fact that $m=2$, there exists $e_0$ and $e_1$ so that the operator $\tilde{\mathcal{L}}^0:=-[1+(k_*\partial_\xi)^2]^2+\mu-3\up^2(\xi;k_*)$, which is the counterpart of the operator $\mathcal{P}$, satisfies
\[
\tilde{\mathcal{L}}^0 e_0=0, \quad \tilde{\mathcal{L}}^0( e_1 + \rmi \xi e_0)=0.
\]
By definition, $\tilde{\mathcal{L}}^0$ is a rescaling of $\mathcal{L}^0$ and thus $e_0$ is the normalized version of $u'_p=k_*\partial_\xi \up$. According to the dependence on parameter $k$ of $\up(\xi; k)$, we readily derive 
\[
\tilde{\mathcal{L}}^0 (\partial_k\up+x\partial_\xi \up)=0,
\]
which, combined with the invertibility of $\tilde{\mathcal{L}}^0$ restricted to the subspace of even, $2\pi$-periodic functions, shows that 
$\partial_k\up+x\partial_\xi \up$ is a rescaling of $e_1 + \rmi \xi e_0$.
As a result, we now conclude that the results in this proposition follows naturally from the self-adjointness of $\mathcal{L}^0$.
\end{Proof}

\paragraph{Spanning the cokernel.}
As a next step, we compute scalar products between
 \[K_1:=\partial_\psi K|_{\psi=0}=(K_{\varphi_0},K_{k_0},K_{\varphi_1},K_{k_1}),\]
 and the elements in the cokernel.  More precisely, we show that $K_{\varphi_0} = K_{k_0} =0$ and that $K_{\varphi_1}$ and $K_{k_1}$ span $\upk$ and $\up'$ in the sense of 
 \begin{equation}\label{e:spandet}
 \det\begin{pmatrix} \langle \up',K_{\varphi_1}\rangle & \langle \upk,K_{\varphi_1}\rangle \\  \langle \up',K_{k_1}\rangle  & \langle \upk,K_{k_1}\rangle \end{pmatrix}\neq0.
 \end{equation}
 where $\langle\cdot, \cdot\rangle$ denotes the standard inner product in $L^2(\R)$.
 
 To start with, a straight forward calculation shows that the total derivative of $K$ is
\begin{equation}\label{e:K1}
 \partial_\psi K|_{\psi =0} =\mathcal{L}^0 (\partial_\xi \up \partial_\psi \varphi |_{\psi =0}+\partial_k \up\partial_\psi \varphi' |_{\psi =0})-\sum_{\pm}\chi_\pm \mathcal{L}^0(\partial_\xi \up \partial_\psi \varphi^\pm |_{\psi =0}+\partial_k \up\partial_\psi (\varphi^\pm)' |_{\psi =0}) 
 \end{equation}
 where $\mathcal{L}_0=L_\mathrm{SH}+F'(\up)$ as defined in \eqref{e:L0} and
 \begin{align*}
 \partial_\psi \varphi &=(-1,x, -\theta,\Theta), & \partial_\psi \varphi' &=(0, 1, -\theta', \theta), \\
 \partial_\psi \varphi^\pm &=(-1,x,\mp 1, \pm x), & \partial_\psi (\varphi^\pm)' &=(0, 1, 0, \pm 1).
 \end{align*}
 We then exploit the fact that $\chi_\pm$ is a partition of unity and $\theta = \chi_+ - \chi_-$ to obtain expressions for each partial derivative in \eqref{e:K1},
\begin{align*}
K_{\varphi_0} &= K_{k_0}=0,\\
K_{\varphi_1} & = [\theta, \mathcal{L}^0] \partial_\xi \up- \mathcal{L}^0( \theta' \partial_k \up),\\
K_{k_1} &= \mathcal{L}^0 \left( \Theta \partial_\xi \up+ \theta \partial_k \up \right) - \theta \mathcal{L}^0( x \partial_\xi \up + \partial_k \up).
\end{align*}
Recalling that $\upk = x \partial_\xi \up + \partial_k \up$, we can further simplify the formula for $K_{k_1}$ into the following form,
\[ K_{k_1} = [\mathcal{L}^0, \theta] \upk +  \mathcal{L}^0 \left( \Theta \partial_\xi \up - \theta x \partial_\xi \up\right).\]
We now proceed to show that \eqref{e:spandet} is true. Noting that $\mathcal{L}^0$ is self-adjoint, $\theta'$ and $\Theta-\theta x$ are compactly supported, 
$\up'= k_* \partial_\xi \up$ and 
\[
[\mathcal{L}^0,w]v=L_\mathrm{SH}(w v) - w L_\mathrm{SH}v=[-\partial_x^4-2\partial_x^2,w]v,
\]
we derive the expressions of projections of $K_{\varphi_1}$ and $K_{k_1}$ on the cokernel,
\begin{align}
\langle \up',K_{\varphi_1}\rangle &= k_*^{-1}\langle \up',[\theta,\mathcal{L}^0]\up'\rangle
=k_*^{-1}\langle \up',[\partial_x^4+2\partial_x^2,\theta]\up'\rangle, \label{e:cok1}\\
\langle \upk,K_{\varphi_1}\rangle &=  k_*^{-1}\langle \upk,[\theta,\mathcal{L}^0]\up'\rangle
=  k_*^{-1}\langle \upk,[\partial_x^4+2\partial_x^2,\theta]\up'\rangle,\label{e:cok2}\\
\langle \up',K_{k_1}\rangle &= \langle \up',[\mathcal{L}^0,\theta]\upk\rangle
=  -\langle \up',[\partial_x^4+2\partial_x^2,\theta]\upk\rangle,\label{e:cok3}\\
\langle \upk,K_{k_1}\rangle &=  \langle \upk,[\mathcal{L}^0,\theta]\upk\rangle
=  -\langle \upk,[\partial_x^4+2\partial_x^2,\theta]\upk\rangle.\label{e:cok4} 
\end{align}
A straightforward computation gives
\begin{equation}\label{e:liebra}
\int_\R u[\partial_x^{2m},w]v\,\rmd x=\int_\R w'\sum_{j=0}^{2m-1}(-1)^{j}u^{(j)}v^{(2m-1-j)}\,\rmd x,
\end{equation}
which has the following two consequences related to \eqref{e:spandet}.
\begin{itemize}
\item[(\rmnum{1})]Applying \eqref{e:liebra} to equation \eqref{e:cok2} and \eqref{e:cok3}, we conclude that the off-diagonal elements in \eqref{e:spandet} coincide, taking the form
\begin{equation}\label{e:offd}
\langle \up',K_{k_1}\rangle=k_*\langle \upk,K_{\varphi_1}\rangle=\int_\R \theta'\left[\sum_{j=0}^{3}(-1)^{j}\upk^{(j)}\up^{(4-j)}+2\sum_{j=0}^{1}(-1)^{j}\upk^{(j)}\up^{(2-j)}\right]\,\rmd x.
\end{equation}
\item[(\rmnum{2})]The expression \eqref{e:liebra} is zero if $u\cdot v\cdot w$ is odd and each of $u, v, w$ is either even or odd. Noting that $ \up'$ and $\theta$ are odd and $\upk$ is even, we conclude that the diagonal elements in \eqref{e:spandet} vanish, that is,
\begin{equation}\label{e:Kd}
\langle \up',K_{\varphi_1}\rangle = \langle \upk,K_{k_1}\rangle=0.
\end{equation}
\end{itemize}
 To further simplify the expression of off-diagonal elements \eqref{e:offd}, we notice that the projections on the cokernel are independent of the choice of $\theta$. More specifically, suppose $\theta_1$ and $\theta_2$ differ by a compactly supported term, $\delta\theta$, we can evaluate the contribution of $\delta\theta$ to our projections,
\[
\int_\R \up' [\mathcal{L}^0,\delta\theta] \upk \,\rmd x= \int_\R \up' \mathcal{L}^0(\delta\theta \upk)-\up' \delta\theta \mathcal{L}^0 \upk\,\rmd x=0.
\]
As a result, the expression in \eqref{e:offd} converges, as $\theta'\to 2\delta_{x_0}$, to 
\begin{align}
\langle \up',K_{k_1}\rangle=k_*\langle \upk,K_{\varphi_1}\rangle=\left.2\left[\sum_{j=0}^{3}(-1)^{j}\upk^{(j)}\up^{(4-j)}+2\sum_{j=0}^{1}(-1)^{j}\upk^{(j)}\up^{(2-j)}\right]\right|_{x=x_0},\label{e:kk}
\end{align}
where $x_0\in\R$ is arbitrary.
Now, using $\upk=\frac{x}{k_*}\up'+\partial_k\up$ and $\up'(0)=\up'(2\pi/k_*)=0$, averaging the constant expression in \eqref{e:kk} over a period $x_0\in [0,2\pi/k_*]$ and integrating by parts, we find,
\begin{equation}\label{e:coker}
\begin{aligned}
\langle \up',K_{k_1}\rangle=k_*\langle \upk,K_{\varphi_1}\rangle
&= \frac{2}{\pi }\int_0^{2\pi/k_*} \left[ k_*\partial_k\left((\up'')^2-(\up')^2\right) + \left( 3(\up'')^2-(\up')^2 \right)  \right] \,\rmd x.
\end{aligned}
\end{equation}
We will see how this expression relates to the effective diffusivity, next, and hence conclude that it does not vanish. As a consequence, $\mathcal{L}^0$ is bounded invertible.

\paragraph{Computing the effective diffusivity.}
We first recall the definition of $L_\mathrm{B}(\sigma)$ from \eqref{e:bl}, and consider the eigenvalue equation
\begin{equation}\label{e:eigen}
L_\mathrm{B}(\sigma)e(\sigma)=\lambda(\sigma)e(\sigma),
\end{equation}
for $\lambda(0)=0$ and $\sigma \sim 0$. Expanding  
\[
L_\mathrm{B}(\sigma)=L_0+L_1\sigma+L_2\sigma^2+\rmO(\sigma^3), \quad e(\sigma)=e_0+e_1\sigma+e_2\sigma^2+\rmO(\sigma^3),  \quad \lambda(\sigma)=\lambda_2\sigma^2 + \rmO(3),
\] 
and setting $e_0=\up'$ and $\langle e_0, e(\sigma)-e_0\rangle_{L^2(0,2\pi/k_*)}=0$, we find explicitly
\[
L_0=-(1+\partial_x^2)^2+\mu -3\up^2(x),\quad L_1=-4\rmi(1+\partial_x^2)\partial_x, \quad L_2=2+6\partial_x^2,
\]
which, plugged in the eigenvalue equation \eqref{e:eigen}, solve 
\[
L_0e_0=0,\quad L_1e_0+L_0e_1=0, \quad L_0e_2+L_1e_1+L_2e_0=\lambda_2e_0.
\] 
Noting $\langle e_1,e_0\rangle_{L^2(0,2\pi/k_*)}=0$, we project the equation for $\lambda_2$ onto $e_1$, that is,
\begin{equation}\label{e:lambda2}
\lambda_2\langle e_0,e_0\rangle_{L^2(0,2\pi/k_*)}=\langle L_1e_1+L_2e_0,e_0\rangle_{L^2(0,2\pi/k_*)}.
\end{equation}
In order to determine $e_1$, we recall Lemma \ref{l:family} and notice that the derivative $\partial_k\up(kx;k)$ at $k=k_*$ satisfies
\[
-4k_*(1+k_*^2\partial_\xi^2)\partial_\xi^2\up +\left(-(1+k_*^2\partial_\xi^2)^2+\mu-3\up^2\right)\partial_k\up=0,
\]
or equivalently, $L_1 e_0+L_0(\rmi k_*\partial_k\up)=0$, which gives 
\[
e_1=\rmi k \partial_k \up.
\]
Inserting the expansion for $L_1$, $L_2$ and $e_1$ into equation \eqref{e:lambda2} gives 
\begin{equation}\label{e:d}
\lambda_2\int_0^{2\pi/k_*}(\up')^2\,\rmd x=-2\int_0^{2\pi/k_*} \left[ k_*\partial_k\left((\up'')^2-(\up')^2\right) + \left( 3(\up'')^2-(\up')^2 \right)  \right] \,\rmd x.
\end{equation}
Therefore, combining \eqref{e:coker} and \eqref{e:d}, we conclude
\begin{equation}\label{e:diffu}
\langle \up',K_{k_1}\rangle=k_*\langle \upk,K_{\varphi_1}\rangle=-\frac{\lambda_2}{\pi}\int_0^{2\pi/k_*} (\up')^2\,\rmd x.
\end{equation}

\begin{Remark}
Notice that a similar reasoning to the proof of Proposition \ref{prop:fredholm} shows that for $\gamma>3/2$ the operators $\mathcal{L}^{\pm,\psi} = L_\textrm{SH} + F'(\up^{\pm, \psi})$, with $\up^{\pm ,\psi}$ as in equation \eqref{e:defup}, are also Fredholm operators from $H^4_{\gamma-2}$ to $L^2_{\gamma}$. Moreover, because the inner products \eqref{e:cok1},\eqref{e:cok2},\eqref{e:cok3}, and \eqref{e:cok4} depend continuously on the parameter $\psi$, the terms $K_{\phi_1}$ and $K_{k_1}$  span the cokernel of these operators as well.
\end{Remark}

\subsection{Invertibility of  $  \mathcal{L}_1^{\psi} $} \label{s:4.2}
The invertibility of $  \mathcal{L}_1^{\psi} $ for $\psi=(0,\varphi_0, 0,0)$ can be derived in a straightforward fashion from the invertibility of $\mathcal{L}_1^0$ due to the simple fact that $ \mathcal{L}_1^{\psi} $ for $\psi=(0,\varphi_0, 0,0)$ is conjugate to $\mathcal{L}_1^0$ via a spatial translation. As a result, we only need to deal with the operator $\mathcal{L}_1^{\psi}$ for $\psi\sim0$. The operators $\mathcal{L}_1^\psi $ are close to $\mathcal{L}^{0}_1$, but the difference is in general not relatively bounded. The difficulty stems from the fact that $ \mathcal{L}^0_1$ ``gains localization'' in certain components, whereas the difference $ \mathcal{L}_1^\psi-\mathcal{L}^0_1$, a bounded multiplication operator, does not affect localization. Therefore, a simple Neumann series perturbation argument will not suffice to establish invertibility of $\mathcal{L}_1^\psi$. We establish somewhat weaker bounds on an inverse of $ \mathcal{L}^\psi_1$ as follows.  First, using the results from subsection \ref{s:4.1} and changing notation in oder to make the distinction between variables and parameters explicit, we write a more complete definition of $\mathcal{L}_1^{\vartheta}$, that is,
\begin{equation}\label{e:lphiinv}
\mathcal{L}^{\vartheta}_1(w,\psi_1):=   -(1+ \partial_x^2)^2 w + \mu w - 3( \up^{\vartheta})^2 w+ K_{\varphi_1} \alpha_0+K_{k_1} \alpha_1 = h
\end{equation}
where $\vartheta = (\vartheta_1,\vartheta_2,\vartheta_3,\vartheta_4)$ denotes the parameter, and  $w,\psi_1=(\alpha_0,\alpha_1)$ are variables. The following proposition then shows the invertibility of this operator and its differentiability with respect to $\vartheta$.

\begin{Proposition}
For $\gamma>3/2$, equation \eqref{e:lphiinv} possesses a solution $(w,\psi_1)$ such that 
\[
\|w\|_{H^4_{\gamma-2}}+|\psi_1|\leq C \|h\|_{L^2_\gamma},
\]
with constant $C$ independent of $\vartheta$, sufficiently small. Moreover, the solution depends continuously on $\vartheta$ in $H^4_{\gamma-2-\delta}$, and is differentiable in $\vartheta$, when considered in spaces with weaker localization,
\[
\|\partial_{\vartheta}w\|_{H^4_{\gamma-3-\delta}}+|\partial_{\vartheta}\psi_1|\leq C \|h\|_{L^2_\gamma}.
\]
\end{Proposition}
\begin{Proof}
For ease of notation we let $m_0 = K_{\varphi_1}, m_1 = K_{k_1}$, and look for solutions to 
\begin{equation}\label{e:Lvarphi} \mathcal{L}_1^\vartheta(w,\psi_1) = \mathcal{L}^\vartheta w + \alpha_0m_0+ \alpha_1 m_1=h,
\end{equation}
where $w\in H^4_{\gamma-2}$, $\alpha_0,\alpha_1\in\R$ are variables, $h \in L^2_{\gamma-2}$, and
\[
\mathcal{L}^{\vartheta}w = -(1+ \partial_x^2)^2 w + \mu w - 3( \up^{\vartheta})^2 w.
\]
We recall as well that $m_0$ and $m_1$ span the cokernel of  
  $ \mathcal{L}^{\pm, \vartheta}= -(1+ \partial_x^2)^2  + \mu  - 3( \up^{\pm, \vartheta})^2,$
   where $\up^{\pm, \vartheta}$ follows the same definition as in equation \eqref{e:defup}. We decompose \eqref{e:Lvarphi} using the partition of unity, $w=w_++w_-$, $h=h_++h_-$, $w_\pm=\chi_\pm w, h_\pm=\chi_\pm h$, and obtain 
\begin{align}
 \mathcal{L}^{ +,\vartheta}w_+ +\sum_{j=0}^1(\alpha_j-\beta_j)m_j+\left(\mathcal{L}^\vartheta-  \mathcal{L}^{-,\vartheta} \right)w_--h_+&=0,\label{e:+}\\
 \mathcal{L}^{-, \vartheta } w_-+\sum_{j=0}^1 \beta_j m_j+\left(\mathcal{L}^\vartheta- \mathcal{L}^{+,\vartheta}\right)w_+-h_-&=0.\label{e:-}
\end{align}
To solve \eqref{e:+} and \eqref{e:-} for $w_\pm$,  $\alpha_j, \beta_j$, $j \in\{0,1\}$, we will consider the cross-coupling terms $\left(\mathcal{L}^\vartheta- \mathcal{L}^{\pm, \vartheta} \right)w_\pm$ as small perturbations. Note that, given $h\in L^2_\gamma$, the system 
\begin{align*}
 \mathcal{L}^{+, \vartheta} w_++\sum(\alpha_j-\beta_j)m_j-h_+&=0\\
 \mathcal{L}^{-, \vartheta}  w_-+\sum \beta_j m_j-h_-&=0,
\end{align*}
possesses a unique solution, $(w_+, w_-, \alpha_1, \alpha_2, \beta_1, \beta_2)$, where $w_-\in H^4_{\gamma-2,\gamma'},w_+\in H^4_{\gamma',\gamma-2}$, with $\gamma'$ arbitrarily large since $h_\pm$ are supported on $\pm x>-1$. Given $|\vartheta|$ small, the cross terms are small, bounded operators when considered on these spaces since, for instance, $\mathrm{supp}(\mathcal{L}^\vartheta-\mathcal{L}^{-,\vartheta} )\subset \R^+$, and $w_-|_{\R^+}\in H^4_{\gamma'}$. This establishes the existence of a bounded inverse, with $w=w_++w_-\in H^4_{\gamma-2}$. It remains to establish the desired smooth dependence of the solution $\underline{w}=(w,\alpha_0,\alpha_1)$  on $\vartheta$.   Writing $\mathcal{L}_1^{\vartheta} \underline{w} = h$ briefly as $\mathcal{L}(\vartheta)(\underline{w}(\vartheta))=h$, we find
\[
\underline{w}(\vartheta+\zeta\varrho)-\underline{w}(\vartheta)=-\mathcal{L}(\vartheta)^{-1}\left(\mathcal{L}(\vartheta+\zeta\varrho)-\mathcal{L}(\vartheta)\right)\underline{w}(\vartheta+\zeta \varrho),
\]
where $0<\zeta\ll 1$, $\vartheta, \varrho\in\R^4$ with $|\varrho|=1$ and $|\vartheta|$ sufficiently small.
Now $\mathcal{L}(\vartheta)^{-1}(\mathcal{L}(\vartheta+\zeta\varrho)-\mathcal{L}(\vartheta) )$ converges to zero when considered as an operator from $H^4_{\gamma-2}\to H^4_{\gamma-2-\delta}$, for any $\delta>0$, which, using uniform bounds for $\underline{w}(\vartheta+\zeta\varrho)$, establishes continuity. Difference quotients and therefore continuity of partial derivatives can be established in a similar fashion.  Notice however that the dependence of the operator $\mathcal{L}^{\vartheta}$ on the parameter comes from the coefficient 
\[3(\up^{\vartheta})^2= 3[\up( k_*x+ \varphi;k_* + \varphi')]^2, \]
via
\[ \varphi(x) = \vartheta_1 x + \vartheta_2 + \vartheta_3\Theta(x) - \vartheta_4 \theta(x).\]
Therefore, derivatives of $\underline{w}(\vartheta)$ with respect to $\vartheta_j$, $j=1,3$ induce linear growth and involve loss of one degree of localization.
\end{Proof}

\subsection{Reduced equations and expansions}\label{s:4.3}
In order to obtain approximations for the variables $(w,\varphi_1, k_1)$, we assume expansions of the form
\begin{align*}
w&=  w_1(\varphi_0,k_0) \varepsilon +  \rmO( \varepsilon^2),\\
\varphi_1& = M_\varphi(\varphi_0,k_0) \varepsilon  +\rmO(\varepsilon^2),\\
k_1 & = M_k(\varphi_0,k_0) \varepsilon + \rmO(\varepsilon^2),
\end{align*}
and we observe that the first order approximations of $(w_1,M_\varphi,M_k)$ satisfy the following equation
\[ \mathcal{L}^0 w_1 + K_{\varphi_1} M_{\varphi}+ K_{k_1} M_k + G_1=0,\]
where by Remark \ref{r:shift} we have that
\[ G_1 = g(x'+\frac{\varphi_0}{k_*}, \up((k_*+k_0)x'; k_*+k_0)).\]
We then proceed to use Lyapunov-Schmidt reduction and obtain the following reduced equations by projecting on the cokernel of $\mathcal{L}^0$,
\begin{align*}
0&= \langle \upk , K_{\varphi_1} \rangle M_\varphi + \langle \upk , G_1\rangle\\
0&= \langle \up' , K_{k_1}\rangle M_k + \langle \up' , G_1\rangle,
\end{align*}
where the variables $M_\varphi$ and $M_k$ depend on $k_0$ and $\varphi_0$. Then, combining these results with  \eqref{e:diffu} and \eqref{e:kk}, and in the particular case of $k_0=0$, we obtain formulas for $M_\varphi(\varphi_0,0)$ and $M_k(\varphi_0,0)$, that is,
\begin{align*}
M_{\varphi}(\varphi_0,0) &= \frac{\pi k_* \displaystyle \int_\R  g(x'+\frac{ \varphi_0}{k_*}, \up)  \upk \,\rmd x'}{\lambda_2 \int_0^{2\pi/k_*} (\up')^2 \,\rmd x}, \\[2ex]
M_{k}(\varphi_0,0)&= \frac{\pi \displaystyle \int_\R g(x'+\frac{ \varphi_0}{k_*}, \up)  \up'  \,\rmd x'}{\lambda_2 \int_0^{2\pi/k_*} (u'_p)^2 \,\rmd x} .
\end{align*}
 It is useful to consider again the change of variables $x' = x - \frac{\varphi_0}{k_*}$, and write
\[ \int_\R g(x'+\frac{ \varphi_0}{k_*}, \up)  \up' \,\rmd x'  = \int_\R  g(x, \up(k_*x - \varphi_0; k_*))  \up' (k_*x - \varphi_0; k_*)\,\rmd x,\]
which, in the case of $g = \partial_u H(x,u)$ for some function $H$, implies that 
\[ \dashint M_k\,\rmd \varphi_0:=\frac{1}{2\pi}\int_0^{2\pi}M_k(\varphi_0,0)\,\rmd \varphi_0=0.\]

\section{Discussion} \label{s:5}

In this paper, we developed a functional-analytic framework for perturbation theory in the presence of essential spectrum, induced by non-compact translation symmetry. The key ingredient are algebraically weighted spaces, including loss of localization by the inverse according to the spatial multiplicity of the essential spectrum. We restricted to ``simple'' branches of essential spectrum for notational simplicity but the methods generalize to more complicated situations. The framework included problems on infinite lattices and cylinders. A crucial assumption is that there is precisely one unbounded direction. 

We showed how such results can be used to study defects, here impurities, in striped phases. The framework of algebraically localized spaces here allows for algebraic decay of impurities. One naturally encounters negative Fredholm indices in the linearization, which one compensates for by adjusting parameters in the far field. In fact, the spatial multiplicity is related in a direct way to the fact that periodic patterns come in two-parameter families. Technically, the decomposition into core deformations (algebraically localized functions) and far field deformations (wavenumber and phase corrections) can be employed in a variety of different contexts. In particular, our approach lays the basis for the continuation of localized deformations such as defects in parameters using more classical algorithms of numerical continuation \cite{lloyd,morrissey}. 

We emphasize that our results do not depend on the particular equation, studied, as long as one is able to determine the existence of periodic patterns and establish properties of the linearization. It is worth noting that both, existence and stability properties, can be established in very reliable ways solving simple periodic boundary-value problems. In particular, one can treat reaction-diffusion systems without much adaptation. Technically more interesting would be systems with conserved quantities such as Cahn-Hilliard, Phase-Field, or DiBlock Copolymer models, since mass conservation induces an additional multiplicity in the essential spectrum, thus violating Hypothesis \ref{h:3.2} on simple kernels of $L(0)$. One could also study problems in channels or infinite cylinders, in particular deformations of hexagonal spot arrays with periodicity of inhomogeneities in one direction. 

There are at least two alternative approaches. First, one could work in exponentially weighted spaces, resorting to stronger assumptions on the inhomogeneity. Fredholm properties of differential operators on the real line in exponentially weighted spaces are well known \cite{palmer,ssmorse} and have been used in the context of perturbation and bifurcation theory in the presence of essential spectrum \cite{ssmorse,goharma}.

In a similar vein, one could cast the existence problem as a non-autonomous differential equation in space $x$, and use dynamical systems tools to investigate the effect of inhomogeneities. From this point of view, the periodic patterns form a two-dimensional normally hyperbolic manifold of equilibria. One can then readily calculate the effect of inhomogeneities on the periodic flow on this center manifold, using traditional methods of averaging.

%
%
%
%
%
%
%
%

A major drawback of these more subtle methods is the reliance on a phase space and exponential behavior in normal directions. In particular, there is no clear path towards perturbation of two-dimensional patterns. Algebraic weights, however, allow for finite-dimensional reductions in the presence of essential spectrum also in higher dimensions \cite{jara1,jara2}. 

\section{Appendix}

\subsection{Fredholm properties of pseudo-derivatives $[D(-\rmi\partial_x)]^{-\ell}$}
\label{ss:a1}

In this section we prove a more general version of Proposition \ref{p:regdrv}. More specifically, for any $\ell\in\Z^+$, 
$p \in (1, \infty)$ and $\gamma_\pm\in\R$, we define the regularized derivative,
\begin{equation}\label{e:rell}
\begin{matrix}
[D(-\rmi\partial_x)]^\ell: & \mathcal{D}([D(-\rmi\partial_x)]^\ell)\subset L^p_{\gamma_--\ell, \gamma_+-\ell} & \longrightarrow & L^p_{\gamma_-, \gamma_+} \\
& u & \longmapsto & \partial_x^\ell ( 1 + \partial_x)^{-\ell}u,
\end{matrix}
\end{equation}
with its domain $\mathcal{D}([D(-\rmi\partial_x)]^\ell)=\{u\in L^p_{\gamma_--\ell, \gamma_+-\ell}\mid (1+\partial_x)^{-\ell}u\in M^{\ell,p}_{\gamma_--\ell, \gamma_+-\ell}\}$. 
From Lemma \ref{p:1ppx} it is straightforward to see that $\mathcal{D}([D(-\rmi\partial_x)]^\ell)$ is a Banach space under the norm
\[
\| u \|:=\| u \|_{L^p_{\gamma_--\ell, \gamma_+-\ell}}+ \| (1+\partial_x)^{-\ell}u\|_{M^{\ell,p}_{\gamma_--\ell, \gamma_+-\ell}}.
\]
Moreover, the Fredholm properties of the bounded operator $[D(-\rmi\partial_x)]^\ell: \mathcal{D}([D(-\rmi\partial_x)]^\ell) \to L^p_{\gamma_-, \gamma_+}$ are summarized in the following proposition. 
\begin{Proposition}\label{p:regdrv}
For $\gamma_\pm\in \mathbb{R} / \{ 1-1/p, 2-1/p, \cdots, \ell-1/p\}$, the regularized derivative $[D(-\rmi\partial_x)]^\ell$ as defined in \eqref{e:rell} is Fredholm. Moreover, the operator $[D(-\rmi\partial_x)]^\ell$ satisfies the following conditions.
\begin{itemize}
 \item If $\gamma_{max}\in I_0:=(-\infty,1-1/p)$, the operator $[D(-\rmi\partial_x)]^\ell$ is onto with its kernel equal to $\mathbb{P}_{\ell}(\R)$.
 \item If $\gamma_{min}\in I_\ell:=(\ell-1/p, \infty)$, the operator $[D(-\rmi\partial_x)]^\ell$ is one-to-one with its cokernel equal to $\mathbb{P}_{\ell}(\R)$.
 \item If $\gamma_{min}\in I_i$, $ \gamma_{max}\in I_j$ with $I_k:=(k-1/p, k+1/p)$ for $0< k\in \Z < \ell$, the kernel and cokernel of the operator $[D(-\rmi\partial_x)]^\ell$ are respectively spanned by $\mathbb{P}_{\ell-j}(\R)$ and $\mathbb{P}_{i}(\R)$.
\end{itemize}
On the other hand, the range of the operator $[D(-\rmi\partial_x)]^\ell$ is not closed if $\gamma_-, \gamma_+\in\{1-1/p, 2-1/p,...,\ell-1/p\}$. 
\end{Proposition}

We will only prove the result in the isotropic case, that is for $\gamma_-= \gamma_+ = \gamma$, since  the proof for the anisotropic case follows the same arguments with straightforward modifications. We start by showing in Lemma \ref{p:1ppx} that the operator $(1 \pm \partial_x) : W^{\ell,p}_{\gamma} \rightarrow W^{\ell-1,p}_{\gamma}$ is an isomorphism and then establish the Fredholm properties of  $\partial_x^\ell:M^{k+\ell,p}_{\gamma-\ell} \rightarrow M^{k,p}_{\gamma}$ in Lemma \ref{p:px}. By combining these two results one arrives at Proposition \ref{p:regdrv}.
\begin{Lemma}
\label{p:1ppx}
Given $\ell\in \Z^+$, $p\in(1,\infty)$, $\gamma\in\R$, the operator $1\pm \partial_x: W^{\ell,p}_{\gamma} \longrightarrow W^{\ell-1,p}_{\gamma}$ is an isomorphism.
\end{Lemma}

\begin{Proof}
 We have the following commutative diagram 
\[
\begin{matrix}
W^{\ell,p}_\gamma &  & \overset{1\pm\partial_x}{\longrightarrow}   &      &W^{\ell-1,p}_\gamma        \\[2mm]
\lfloor x \rfloor^\gamma\downarrow    &     &                         &          &\lfloor x \rfloor^\gamma\downarrow\\
W^{\ell,p}  &     & \overset{\mathcal{M}_\pm}{\longrightarrow}        &       &W^{\ell-1,p}.
\end{matrix}
\]
As a result, we have $(\mathcal{M}_\pm u)(x)= \lfloor x \rfloor^\gamma (1\pm\partial_x)(\lfloor x \rfloor^{-\gamma}u(x))= (1\pm\partial_x)u(x)-\gamma x\lfloor x \rfloor^{-2}u(x)$, that is,
according to the Kondrachov embedding theorem, the operator $\mathcal{M}_\pm$ is equal to a compact perturbation of the  invertible operator $(1\pm \partial_x): W^{\ell,p} \rightarrow W^{\ell-1,p}$. Noting that $\ker{\mathcal{M}_\pm}=\{0\}$, we conclude that $\mathcal{M}_\pm$ is invertible.
\end{Proof}

To obtain the Fredholm properties of  $\partial_x^\ell$, we first generalize the canonical definition of $\partial_x: M^{k+1,p}_{\gamma-1} \rightarrow M^{k,p}_{\gamma}$ where $k\geq 0$ to the $k<0$ regime: given $k\in\Z^-$, the operator  $\partial_x :M^{k+1,p}_{\gamma-1} \rightarrow M^{k,p}_{\gamma}$ is defined as 
\begin{equation}\label{e:gpx}
 \partial_x u (v) = - \llangle u, \partial_x v \rrangle, \quad \forall u \in M^{k+1,p}_{\gamma-1}, v \in M^{-k,q}_{-\gamma},
\end{equation} 
where $1/p+1/q=1$. 
\begin{Remark}
The generalized operator $\partial_x :L^p_{\gamma-1} \rightarrow M^{-1,p}_{\gamma}$ is an extension of the canonical operator $\partial_x: M^{1,p}_{\gamma-1} \rightarrow L^p_{\gamma}$ in the sense that $\partial_xu(v) = \llangle \partial_x u , v \rrangle$, for any $ u  \in M^{1,p}_{\gamma-1}$ and $ v \in M^{1,q}_{-\gamma}$.
\end{Remark}

For this generalized operator, we have the following lemma whose proof will occupy the rest of this section.
\begin{Lemma}\label{p:px}
Given $k\in\mathbb{Z}$, $\ell\in\mathbb{Z}^+$, $p\in(1,\infty)$, and $\gamma\in \mathbb{R}\setminus\{1-1/p, 2-1/p,...,\ell-1/p\}$, the operator
\begin{equation}\label{e:partialxm}
\partial_x^\ell: M^{k+\ell,p}_{\gamma-\ell}\longrightarrow M^{k,p}_{\gamma},
\end{equation}
is Fredholm. Moreover, 
\begin{itemize}
 \item if $\gamma<1-1/p$, the operator \eqref{e:partialxm} is onto with its kernel equal to $\mathbb{P}_\ell(\R)$;
 \item if $\gamma>\ell-1/p$, the operator \eqref{e:partialxm} is one-to-one with its cokernel equal to $\mathbb{P}_\ell(\R)$;
 \item if $j-1/p<\gamma<j+1-1/p$, where $j \in \Z^+\cap[1,\ell-1]$, the kernel  and cokernel of the operator \eqref{e:partialxm} 
 are respectively $\mathbb{P}_{\ell-j}(\R)$ and $\mathbb{P}_j(\R)$.
\end{itemize}
On the other hand, the operator \eqref{e:partialxm} does not have a closed range if $\gamma\in\{1-1/p, 2-1/p,...,\ell-1/p\}$. 
\end{Lemma}

We focus on the proof of the two primary cases when $\ell=1$ and $k=0,-1$, which can be readily generalized to the case when $\ell=1$ and $k=n, -n-1$ for $n\in\Z^+$, and then the case $\ell>1$. The proof is given in various steps written as lemmas. 
We first establish Fredholm properties of the operator  $\partial_x: M^{1,p}_{\gamma-1} \rightarrow L^p_{\gamma}$ when $ \gamma>1-1/p$ in Lemma \ref{l:k1}.  We then establish
Fredholm properties of the operator  $\partial_x: L^p_{\gamma-1} \rightarrow M^{-1,p}_{\gamma}$ when $\gamma\neq 1-1/p$ in Lemma \ref{l:k0}-\ref{l:k01}, where Fredholm properties of the operator $\partial_x: M^{1,p}_{\gamma-1} \rightarrow L^p_{\gamma}$ when $ \gamma<1-1/p$ follow.
Finally, we show in Lemma \ref{l:1-1/p} that for $\gamma =1-1/p$ both operators do not have closed range.
\begin{Lemma}\label{l:k1}
Given $p \in (1,\infty)$ and $\gamma > 1-1/p$, the operator, $ \partial_x: M^{1,p}_{\gamma-1} \rightarrow L^p_{\gamma}$, is Fredholm and one-to-one with its cokernel spanned by $\mathbb{P}_1(\R)$.
\end{Lemma}

\begin{Remark}
We can readily apply the techniques from the following proof to show that, given $p \in (1,\infty)$ and $[\gamma_+ - (1-1/p)][\gamma_--(1-1/p)]<0$, the operator, $ \partial_x: M^{1,p}_{\gamma_--1,\gamma_+-1} \rightarrow L^p_{\gamma_-,\gamma_+}$, is bounded and invertible.
\end{Remark}

\begin{Proof}
Given $\gamma>1-1/p$, we write 
\[
L^p_{\gamma,\perp}:=\{f\in L^p_\gamma\mid \int_\R f=0\},
\]
which is closed in $L^p_\gamma$ since $1$ is a bounded linear functional on $L^p_\gamma$.
It is not hard to see that, for any $u\in M^{1,p}_{\gamma-1}$, its derivative $\partial_x u\in L^1$.
We then consider $v(x):=\int_{\infty}^x\partial_x u(y)\rmd y$ and take $C_1=\lim_{x\to-\infty}v(x)$. It is clear that there exists some $C_2\in \mathbb{R}$ such that
$u(x)-v(x)=C_2,$
which leads to 
$$\lim_{x\to\infty}u(x)=C_2,\quad \lim_{x\to-\infty}u(x)=C_2+C_1.$$
The fact that $u\in L^p_{\gamma-1}$ implies that if the $\lim_{x\to \pm\infty}u(x)$ exists, it must be zero. Thus, we have
$C_1=C_2=0$, that is, $\int_{\mathbb{R}}\partial_x u \rmd x=0$,
and consequently 
\[
\Rg (\partial_x )\subseteq L^p_{\gamma,\perp}.
\]

We now claim that the inverse of $\partial_x$ can be defined as
\begin{equation}\label{e:inv}
\begin{aligned}
\partial_x^{-1}:&\quad L^p_{\gamma,\perp} &\longrightarrow & \quad M^{1,p}_{\gamma-1}\\
\quad &\quad\quad f & \longmapsto & \quad \int_{\infty}^x f(y)\rmd y.
\end{aligned}
\end{equation}
The fact that  $\partial_x^{-1}$ is well defined reduces to verifying that $u(x)=\int_{\infty}^x f(y)\rmd y\in L^p_{\gamma-1}$. To do that, we let $\tilde{\gamma}:=\gamma-(1-1/p)>0$ and split $\R$ into three intervals, that is, $\R=(-\infty,-1)\cup [-1,1]\cup (1,\infty)$.
First, it is not hard to see that
\begin{equation}\label{e:ineq1}
\|u(x)\|_{L^p_{\tilde{\gamma}-1/p}([-1,1])}\leq C(\gamma,p)\max_{|x|\leq 1}|u(x)|
\leq C(\gamma,p)\|f\|_{L^1(\R)}\leq C(\gamma,p)\|f\|_{L^p_\gamma(\R)},
\end{equation}
where $C(\gamma)$ is a constant varying with $\gamma$ and $p$.
For the interval $(1,\infty)$, we use a logarithmic scaling, that is,
\[
\tau:=\ln (x),\quad w(\tau):=\rme^{\tilde{\gamma}\tau}u(\rme^\tau),\quad g(\tau):=\rme^{(\tilde{\gamma}+1)\tau}f(\rme^\tau),
\]
so that the ODE $w_\tau-\tilde{\gamma}w=g$ admits a solution $w(\tau)=\int_\infty^\tau \rme^{\tilde{\gamma}(\tau-s)}g(s)\rmd s$.
Applying Young's inequality  to the above integral equation, we obtain 
\begin{equation}\label{e:ineq2}
\sqrt{2}^{(1/p-\tilde{\gamma})}\|u(x)\|_{L^p_{\tilde{\gamma}-1/p}((1,\infty))}\leq \|w(\tau)\|_{L^p((0,\infty))}\leq 
\frac{1}{\tilde{\gamma}}\|g(\tau)\|_{L^p((0,\infty))} \leq \frac{1}{\tilde{\gamma}}\|f(x)\|_{L^p_{\tilde{\gamma}+1-1/p}((1,\infty))}.
\end{equation}
For the interval $(-\infty,1)$, a similar argument can be applied and leads to the inequality,
\begin{equation}\label{e:ineq3}
\|u(x)\|_{L^p_{\tilde{\gamma}-1/p}((-\infty,-1))}\leq C(\gamma, p)\|f(x)\|_{L^p_{\tilde{\gamma}+1-1/p}((-\infty,-1))}.
\end{equation}
Combining the inequalities (\ref{e:ineq1})--(\ref{e:ineq3}), we conclude that the operator \eqref{e:inv} is well defined and we have
$$\|\partial_x^{-1}f\|_{M^{1,p}_{\gamma-1}}=\|u\|_{L^p_{\gamma-1}}+\|f\|_{L^p_{\gamma}}\leq C(\gamma)\|f\|_{L^p_{\gamma}},$$
which implies that $\partial_x^{-1}$ is also a bounded linear operator. 
\end{Proof}

\begin{Lemma}\label{l:k0}
Given $p \in (1, \infty)$, we have that,
\begin{itemize}
\item for  $\gamma >1-1/p$ , the operator $\partial_x : L^p_{\gamma-1} \rightarrow M^{-1,p}_{\gamma}$ is one-to-one;
\item for $\gamma < 1-1/p$, the operator $\partial_x : L^p_{\gamma-1} \rightarrow M^{-1,p}_{\gamma}$ is Fredholm, onto with its kernel equal to $\mathbb{P}_1(\R)$.
\end{itemize}
\end{Lemma}
\begin{Proof}
For $\gamma> 1-1/p$, consider $u \in L^p_{\gamma-1}$ with $\partial_xu =0$. We let $\{ u_n\}_{n \in \mathbb{N}} \subset C^{\infty}_0$ such that $u_n \rightarrow u$ in $L^p_{\gamma-1}$ and then have that, for any $v\in M^{1,q}_{-\gamma}$,
\[ \partial_x u (v) =- \llangle u, \partial_x v \rrangle = \lim_{n \rightarrow \infty} \llangle \partial_x u_n, v \rrangle =0,\]
which implies $\partial_x u_n \rightarrow 0$ in $L^p_{\gamma}$. We therefore have $u=0$, proving the first statement of the lemma. 

For $\gamma< 1-1/p$, the operator $\partial_x: M^{1,q}_{-\gamma} \rightarrow L^q_{1-\gamma}$, according to Lemma \ref{l:k1}, is a Fredholm operator with index $-1$ and cokernel  equal to $\mathbb{P}_1(\R)$. Therefore, the operator $\partial_x: L^p_{\gamma-1} \rightarrow M^{-1,p}_{\gamma}$, as the adjoint operator of $\partial_x: M^{1,q}_{-\gamma} \rightarrow L^q_{1-\gamma}$ with an extra negative sign, is Fredholm with index $1$ and kernel equal to $\mathbb{P}_1(\R)$.
\end{Proof}

\begin{Lemma}\label{l:k01}
Given $p \in (1, \infty)$, we have
\begin{itemize}
\item for $\gamma< 1-1/p$, the Fredholm operator $\partial_x: M^{1,p}_{\gamma-1} \rightarrow L^p_{\gamma}$ is onto with its kernel equal to $ \mathbb{P}_1(\R)$.
\item for $\gamma>1-1/p$, the Fredholm operator $\partial_x: L^p_{\gamma-1} \rightarrow M^{-1,p}_{\gamma}$ is one-to-one with its cokernel equal to $ \mathbb{P}_1(\R)$.
\end{itemize}
\end{Lemma}
\begin{Proof}
To prove the lemma we just need to show that each operator has a closed range. We restrict our attention to the first operator, the second being analogous. By way of contradiction, suppose that $\partial_x: M^{1,p}_{\gamma-1} \rightarrow L^p_{\gamma}$ does not have a closed range for $\gamma<1-1/q$, then there exists a sequence $\{ u_n \}_{n \in \mathbb{N}} \subset M^{1,p}_{\gamma-1}$ such that dist$(u_n , \mathbb{P}_1(\R)) =1$ and $\| \partial_x u_n \|_{L^p_{\gamma}} \rightarrow 0$. The norm inequality $\| \partial_x u_n \|_{M^{-1,p}_{\gamma}} \leq \|\partial_x u_n \|_{L^p_{\gamma}}$, together with the fact that the operator $\partial_x: L^p_{\gamma-1} \rightarrow M^{-1,p}_{\gamma}$ has closed range show that we can find a subsequence $\{ v_n\} \subset \ker(\partial_x) \subset M^{1,p}_{\gamma-1}$ such that $\|u_n -v_n\|_{L^p_{\gamma-1}} \rightarrow 0$. Therefore, we have
\[
\|u_n -v_n\|_{M^{1,p}_{\gamma-1}} \leq \|u_n -v_n\|_{L^p_{\gamma-1}} + \|\partial_x u_n - \partial_x v_n\|_{L^p_{\gamma}}\rightarrow 0, \text{ as }n\to\infty,
\]
that is, $\mathrm{dist}(u_n , \mathbb{P}_1(\R)) \rightarrow 0$, which is a contradiction and concludes the proof.
\end{Proof}

\begin{Lemma}\label{l:1-1/p}
Given $p \in (1,\infty)$ and $\gamma=1-1/p$, the operators $\partial_x : M^{1,p}_{\gamma-1} \rightarrow L^p_{\gamma}$ 
and $\partial_x: L^p_{\gamma-1} \rightarrow M^{-1,p}_{\gamma}$ do not have closed range.
\end{Lemma}
\begin{Proof}
 Let $\phi \in C^{\infty}_0$ with $0 \leq \phi \leq 1$ and $\mathrm{supp}(\phi)=[-1,1]$.
Let $u_n(x) = \phi(x/n)$, then $\{\partial_x u_n\}_{n\in\Z^+}$ is a bounded sequence in $M^{-1,p}_{\gamma}$ (also, in $L^p_\gamma$). However, if $\gamma=1-1/p$ the sequence $\{ u_n \}_{n \in \mathbb{N}}$ is unbounded in $L^p_{\gamma-1}$ (also, in $W^{1,p}_{\gamma-1}$). Therefore, both operators do not have closed range.
\end{Proof}

\subsection{Fredholm properties of operators $\delta_+^{\ell-i}\delta_-^i$}
\label{ss:a2}

\begin{Proposition}\label{p:delta}
Given $k\in\mathbb{Z}$, $\ell\in\mathbb{Z}^+$, $p\in(1,\infty)$, and $\gamma\in \mathbb{R}\setminus\{1-1/p, 2-1/p,...,\ell-1/p\}$, the operator
\begin{equation}\label{e:deltafred}
\delta_+^{\ell-i}\delta_-^i: \mathcal{M}^{k+\ell,p}_{\gamma-\ell}\longrightarrow \mathcal{M}^{k,p}_{\gamma},
\end{equation}
is Fredholm for $i\in[0,\ell]\cap\Z$. Moreover, 
\begin{itemize}
 \item if $\gamma<1-1/p$, the operator in \eqref{e:deltafred} is onto with its kernel equal to $\mathbb{P}_\ell(\Z)$;
 \item if $\gamma>\ell-1/p$, the operator in \eqref{e:deltafred} is one-to-one with its cokernel equal to $\mathbb{P}_\ell(\Z)$;
 \item if $j-1/p<\gamma<j+1-1/p$, where $j \in \Z^+\cap[1,\ell-1]$, the kernel  and cokernel of the operator in \eqref{e:deltafred} 
 are respectively $\mathbb{P}_{\ell-j}(\Z)$ and $\mathbb{P}_j(\Z)$.
\end{itemize}
On the other hand, the operator in \eqref{e:deltafred} does not have a closed range if $\gamma\in\{1-1/p, 2-1/p,...,\ell-1/p\}$. 
\end{Proposition}

The proof of Proposition \ref{p:delta} is essentially the same as in the continuous case, that is, the proof of Lemma \ref{p:px}. The main technical difference lies in the proof of the the discrete version of Lemma \ref{l:k1}, which we shall establish now. 
\begin{Lemma}
For $\gamma>1-1/p$ and $p\in[1,\infty]$, discrete derivative operators, $\delta_\pm: \mathscr{M}^{1,p}_{\gamma-1} \mapsto \ell^p_\gamma$,
are one-to-one Fredholm operators with both cokernels spanned by $\p_1(\Z)$.
\end{Lemma}

\begin{Proof}
It is straightforward to see that $\delta_\pm$ are isomorphic and we only need to prove the results for $\delta_+$.
Just like the continuous, the essential part is to prove that 
\begin{equation*}
\begin{matrix}
 \delta_+^{-1}: & \ell^p_{\gamma,\perp} & \longrightarrow & \ell^p_{\gamma-1} \\
 &  \{b_j\}_{j\in\Z} & \longmapsto & \{-\sum_{i=j}^\infty b_{i}\}_{j\in\Z},
\end{matrix}
\end{equation*}
where $\ell^p_{\gamma,\perp}=\{\{b_j\}_{j\in\Z}\in\ell^p_\gamma\mid \sum_{j\in\Z}b_j=0\}$, is the bounded inverse of $\delta_+$.
To do that, we instead consider the following operator
\begin{equation*}
\begin{matrix}
\widetilde{\delta}_+^{-1}: & \ell^p_{\gamma,\perp}(\N) & \longrightarrow & \ell^p_{\gamma-1}(\N) \\
 &  \{b_j\}_{j\in\N} & \longmapsto & \{-\sum_{i=j}^\infty |b_{i}|\}_{j\in\N},
\end{matrix}
\end{equation*}
We denote $a_j=-\sum_{i=j}^\infty b_{i}$ for all $j\in \Z$ and $\widetilde{a}_j=-\sum_{i=j}^\infty |b_{i}|$ for all $j\in\N$. 
It is then not hard to conclude that
\begin{itemize}
\item $a_{j+1}-a_j=b_j$, for all $j\in\Z$;
\item $\widetilde{a}_{j+1}-\widetilde{a}_j=|b_j|$, for all $j\in\N$;
\item $\{\widetilde{a}_j\}_{j\in \N}$ is an increasing sequence with non-negative entries;
\item $|\widetilde{a}_j|\geq |a_j|$, for all $j\in\N$.
\end{itemize}
For any $\widetilde{\gamma}>0$ and $j\in\N$, we introduce
\begin{equation*}
A_j=2^{j\widetilde{\gamma}}\widetilde{a}_{2^j},\quad B_j=2^{j\widetilde{\gamma}}\sum_{i=2^j}^{2^{j+1}-1}|b_j|,
\end{equation*}
and have $2^{-\widetilde{\gamma}}A_{j+1}-A_{j}=B_j$, or equivalently, $A_j=-\sum_{i=j}^\infty 2^{(j-i)\widetilde{\gamma}}B_i$, which, according to Young's inequality, leads to that
\begin{equation}\label{e:ABI}
 \|\{A_j\}_{j\in\N}\|_{\ell^p(\N)}\leq 
\|\{2^{-\widetilde{\gamma}j}\}_{j\in\N}\|_{\ell^1}\|\{B_j\}_{j\in\N}\|_{\ell^p(\N)}\leq 
\frac{2^{\widetilde{\gamma}}}{2^{\widetilde{\gamma}}-1}\|\{B_j\}_{j\in\N}\|_{\ell^p(\N)}.
\end{equation}
Moreover, on the one hand, we have
\begin{equation}\label{e:AI}
\begin{aligned}
\|\{A_j\}_{j\in\N}\|_{\ell^p(\N)}^p 
= \sum_{j=0}^\infty 2^{\widetilde{\gamma}pj-j}\left(2^j|\widetilde{a}_{2^j}|^p\right)
&\geq \sum_{j=0}^\infty 2^{(\widetilde{\gamma}p-1)j}\left(\sum_{i=2^j}^{2^{j+1}-1}|\widetilde{a}_{i}|^p\right)\\
&\geq \min\{4^{1-\widetilde{\gamma}p},1\}\sum_{j=0}^\infty \left(\sum_{i=2^j}^{2^{j+1}-1}\lfloor i\rfloor^{\widetilde{\gamma}p-1}|\widetilde{a}_{i}|^p\right)\\
&=\min\{4^{1-\widetilde{\gamma}p},1\}\|\{\widetilde{a}_j\}_{j\in\Z^+}\|_{\ell^p_{\widetilde{\gamma}-1/p}(\Z^+)}^p\\
&\geq\min\{4^{1-\widetilde{\gamma}p},1\}\|\{a_j\}_{j\in\Z^+}\|_{\ell^p_{\widetilde{\gamma}-1/p}(\Z^+)}^p.
\end{aligned}
\end{equation}
On the other hand, we have
\begin{equation}\label{e:BI}
\begin{aligned}
\|\{B_j\}_{j\in\N}\|_{\ell^p(\N)}^p =\sum_{j=0}^\infty 2^{(\widetilde{\gamma}+1)pj}\left(\frac{1}{2^j}\sum_{i=2^j}^{2^{j+1}-1}|b_{i}|\right)^p
&\leq \sum_{j=0}^\infty 2^{[(\widetilde{\gamma}+1)p-1]j}\left(\sum_{i=2^j}^{2^{j+1}-1}|b_{i}|^p\right)\\
&\leq \max\{4^{1-(\widetilde{\gamma}+1)p},1\}\sum_{j=0}^\infty \left(\sum_{i=2^j}^{2^{j+1}-1}i^{(\widetilde{\gamma}+1)p-1}|b_{i}|^p\right)\\
&=\max\{4^{1-(\widetilde{\gamma}+1)p},1\}\|\{b_j\}_{j\in\Z^+}\|_{\ell^p_{\widetilde{\gamma}+1-1/p}(\Z^+)}^p.
\end{aligned}
\end{equation}
Combining these inequalities \eqref{e:ABI}, \eqref{e:AI} and \eqref{e:BI}, we conclude that, there exists $C(\tilde{\gamma},p)>0$ so that
\begin{equation*}
\|\{a_j\}_{j\in\Z^+}\|_{\ell^p_{\widetilde{\gamma}-1/p}(\Z^+)}\leq C(\tilde{\gamma},p)
\|\{b_j\}_{j\in\Z^+}\|_{\ell^p_{\widetilde{\gamma}+1-1/p}(\Z^+)}\leq C(\tilde{\gamma},p)
\|\{b_j\}_{j\in\Z}\|_{\ell^p_{\widetilde{\gamma}+1-1/p}(\Z)}.
\end{equation*}
By shifting and letting $j\rightarrow -j$, we can also show that 
\begin{equation*}
\|\{a_j\}_{j\in\Z^-\cup\{0\}}\|_{\ell^p_{\widetilde{\gamma}-1/p}(\Z^-\cup\{0\})}\leq C(\tilde{\gamma},p)
\|\{b_j\}_{j\in\Z}\|_{\ell^p_{\widetilde{\gamma}+1-1/p}(\Z)}.
\end{equation*}
In conclusion, letting $\tilde{\gamma}=\gamma-1-1/p>0$, there exists $C(\gamma,p)>0$ such that
\begin{equation*}
\|\{a_j\}_{j\in\Z}\|_{\ell^p_{\gamma-1}}\leq C(\gamma, p) \|\{b_j\}_{j\in\Z}\|_{\ell^p_{\gamma}},
\end{equation*}
which concludes the proof.
\end{Proof}


\end{document}